\documentclass[a4paper,11pt]{article}

\usepackage[T1]{fontenc}
\usepackage[latin1]{inputenc}
\usepackage{amsmath,amssymb}
\usepackage{theorem}

\newtheorem{thm}{Theorem}[section]
\newtheorem{rem}[thm]{Remark}
\newtheorem{prop}[thm]{Proposition}
\newtheorem{lemma}[thm]{Lemma}
\newtheorem{cor}[thm]{Corollary}
\theoremstyle{break}
\newtheorem{thmbr}[thm]{Theorem}
\newtheorem{lembr}[thm]{Lemma}
\newtheorem{rembr}[thm]{Remarks}

\numberwithin{equation}{section}
\numberwithin{figure}{section}

\newenvironment{proof}{\paragraph{Proof}}
{\hfill$\Box$
\bigskip}

\newcommand{\RR}{\mathbb{R}}

\newcommand{\PP}{\mathbb{P}}
\newcommand{\EE}{\mathbb{E}}

\title{Large deviations for singular and degenerate diffusion models in
  adaptive evolution}

\author{Nicolas Champagnat\thanks{EPI TOSCA, INRIA Sophia Antipolis
    -- M\'editerran\'ee, 2004 route des
Lucioles, BP. 93, 06902 Sophia Antipolis Cedex, France; e-mail:
Nicolas.Champagnat@sophia.inria.fr.}}

\date{\today}

\begin{document}

\maketitle

\begin{abstract}
  In the course of Darwinian evolution of a population, punctualism is
  an important phenomenon whereby long periods of genetic stasis
  alternate with short periods of rapid evolutionary change. This
  paper provides a mathematical interpretation of punctualism as a
  sequence of change of basin of attraction for a diffusion model of
  the theory of adaptive dynamics. Such results rely on large
  deviation estimates for the diffusion process. The main difficulty
  lies in the fact that this diffusion process has degenerate and
  non-Lipschitz diffusion part at isolated points of the space and
  non-continuous drift part at the same points. Nevertheless, we are
  able to prove strong existence and the strong Markov property for
  these diffusions, and to give conditions under which pathwise
  uniqueness holds. Next, we prove a large deviation principle
  involving a rate function which has not the standard form of
  diffusions with small noise, due to the specific singularities of
  the model.
  Finally, this result is used to obtain asymptotic estimates for the
  time needed to exit an attracting domain, and to identify the points
  where this exit is more likely to occur.
\end{abstract}

\noindent {\emph{AMS 2000 subject classification.} Primary 60F10, 92D15;
  secondary 60J70, 60J60.
\medskip
  
}

\section{Introduction}
\label{sec:intro-LDP}

The Darwinian evolution of an asexual population is controlled by
demographic (birth and death) rates, which are typically influenced by
quantitative characters, called phenotypic traits: morphological
traits like body size, physiological traits like the rate of food
intake, life-history traits like the age at maturity. Such traits are
heritable yet not perfectly transmitted from parents to offsprings,
due to mutations of genes involved in their expression. The resulting
variation of traits is then exposed to selection caused by ecological
interactions between individuals competing for limited
resources. Models of evolution of the dominant trait in the space of
phenotypic traits are usually of two types: jump processes (often
called ``adaptive random walks''~\cite{metz-geritz-al-96,drossel-01})
or diffusion processes (\cite{lande-76,hansen-97}). Diffusion models
are usually more suited to finite populations, weak selection, or long
time scales. These models usually involve a so-called ``fitness
function'', which quantifies the selective ability of each possible
phenotypic traits. Such models are also sometimes referred to as
evolution models on a ``fitness landscape'' (an notion going back to
Wright~\cite{wright-32}).

In most cases, the parameters of these models (speed of evolution,
fitness function,\ldots) are based on heuristic
considerations. However, since the early 1990's, adaptive dynamics
theory~\cite{hofbauer-sigmund-90,marrow-law-al-92,metz-nisbet-al-92}
has been developed to give a firm basis to such models, starting from
an individual-based description of the population with explicit
ecological interactions. The combination of ecology and evolution
allowed to obtain evolutionary models on a fitness landscape that
depends on the current state of the population, and which is
explicitly given in terms of individual parameters. The first model
is an adaptive random walk, called the ``trait substitution sequence''
(TSS), first described in~\cite{metz-geritz-al-96} (see
also~\cite{dieckmann-law-96}). The mathematical derivation of this
model from an individual-based model under specific asymptotics has
been done in~\cite{champagnat-06}. In the limit of small mutations,
this stochastic jump process converges to a deterministic ordinary
differential equation called ``canonical equation of adaptive
dynamics''~\cite{dieckmann-law-96,champagnat-ferriere-al-01,champagnat-meleard-08}.
Several diffusion models have also been obtained in this
framework~\cite{champagnat-ferriere-al-06,champagnat-lambert-07},
either as diffusion approximations of the TSS or in the case of weak
selection in finite populations.

One evolutionary pattern that remains poorly understood among
biologists is that of ``punctualism'': the phenomenon of Darwinian
evolution whereby long periods of trait stasis alternate with periods
of global, rapid changes in the trait values of the population, which
can be due to a large mutation or to successive invasions of slightly
disadvantaged mutants in the population~\cite{rand-wilson-93}.
In this paper, we interpret punctualism as phases of quick changes of
basin of attraction for the canonical equation of adaptive dynamics,
separated by long phases where the population state stays near the
evolutionary equilibrium inside the current basin of attraction
(``problem of exit from a domain''~\cite{freidlin-wentzell-84}). The
TSS model is not well-suited to this study because it cannot jump in
the direction of less fitted traits (i.e.\ traits having negative
fitness). However, for punctualism to occur, a sequence of surviving
unfitted mutations must occur. This is possible on long time scales
because of the finiteness of the population. Therefore, we focus in
this work on a diffusion model of adaptive dynamics that generalizes
the one of~\cite{champagnat-ferriere-al-06} (see~\cite{champagnat-04}
for a general derivation of these models), where evolution can proceed
in any direction of trait space. This model is obtained as a diffusion
approximation (in the sense of~\cite[Ch.\:11]{ethier-kurtz-86}) of the
TSS.

This diffusion process on the trait space, assumed to be a
subset of $\RR^d$, is solution to the the following stochastic
differential equation, with coefficients explicitly obtained in terms
of biological parameters (see section~\ref{sec:a-b-btilde}):
\begin{equation}
  \label{eq:SDE-LDP}
  dX^{\varepsilon}_t=(b(X^{\varepsilon}_t)+\varepsilon
  \tilde{b}(X^{\varepsilon}_t))dt+
  \sqrt{\varepsilon}\sigma(X^{\varepsilon}_t)dW_t,
\end{equation}
where $b(x)$ and $\tilde{b}(x)$ are in $\mathbb{R}^d$, $\sigma(x)$ is
a $d\times d$ symmetric positive real matrix, and $\varepsilon>0$ is a
small parameter scaling the size of mutation jumps.

The main difficulty of this model is that the standard regularity
assumptions for stochastic differential equations (SDE) are not
satisfied: the function $b$ is (globally) Lipschitz, but $\tilde{b}$
is discontinuous at isolated points of the trait space, called
\emph{evolutionary singularities}, and $\sigma$ is not globally
Lipschitz, but is only $1/2$-H\"older near the set $\Gamma$ of
evolutionary singularities. Moreover, $b(x)=\tilde{b}(x)=\sigma(x)=0$
for $x\in\Gamma$.

Despite these difficulties, we are able to study the existence, strong
Markov property and (partly) uniqueness for this SDE, to prove a large
deviations principle (LDP) as $\varepsilon\rightarrow 0$, and to study the
problem of diffusion exit from a domain of Freidlin and
Wentzell~\cite{freidlin-wentzell-84}, which is the key question for
punctualism: what are the time and point of exit of $X^\varepsilon$
from an attracting domain?

The original method for proving a LDP for the solution to a SDE with
Lipschitz coefficients was based on discretization and continuous
mapping techniques~\cite{freidlin-wentzell-84,azencott-80} (transfer
of the LDP for Brownian motion---Schilder's theorem---to the LDP for
the diffusion). This technique has been extended to weaker assumptions
on the coefficients (e.g.\ essentially locally-lipschitz
in~\cite{baldi-chaleyat-maurel-88} or for a restricted class of
two-dimensional diffusions in~\cite{korostelev-leonov-93}) or to
reflected diffusions~\cite{doss-priouret-83}. Other techniques were
more recently developed to study LDP for diffusions with irregular
coefficients. The weak convergence approach of Dupuis and
Ellis~\cite{dupuis-ellis-97} is based on a combination of perturbation
approach, discretization and representation formulas. They were in
particular able to obtain upper bounds under very general
assumptions~\cite{dupuis-ellis-al-91} and to obtain the LDP for
diffusions with discontinuous coefficients~\cite{boue-dupuis-al-00}
(see also~\cite{chiang-sheu-00}). Another technique developed by de
Acosta~\cite{deAcosta-00}, is based on an abstract non-convex
formulation of LDP, and allows one to deal with degenerate diffusion
coefficients, but requires Lipschitz coefficients.

However, the existing results dealing with discontinuous coefficients
are of a different nature as the singularity we consider
(in~\cite{chiang-sheu-00,boue-dupuis-al-00}, the drift coefficient is
discontinuous on a hyperplane), and these later methods require either
the coefficients to be Lipschitz, or the diffusion parameter to be
non-degenerate. Another reason why these methods seem not to apply
easily to our situation is that the rate function arising naturally
with these methods does not take into account the singularity of our
model. Actually, the results of~\cite{dupuis-ellis-al-91} can be used
to obtain an large deviation upper bound, but, as appears in
Section~\ref{sec:LDP}, with a non-optimal rate function. For these
reasons, we adapt in this work the original methods based on
discretization and path comparisons, allowing us to finely study the
paths of the diffusion $X^\varepsilon$ near $\Gamma$. Our proof
follows the method of Azencott~\cite{azencott-80} (see
also~\cite{doss-priouret-83}). Interestingly, it also appears that, in
contrast with what is usually observed in large deviations theory (see
e.g.~\cite{dupuis-ellis-al-91}), our upper bound is more difficult to
obtain than the lower bound.

The paper is organized as follows. In Section~\ref{sec:a-b-btilde}, we
describe precisely the model and study the regularity of the
parameters $a=\sigma\sigma^*$, $b$ and $\tilde{b}$. In
Sections~\ref{sec:construction-X} and~\ref{sec:uniqueness-Markov}, we
establish strong existence and the strong Markov property
for~(\ref{eq:SDE-LDP}), by explicitly constructing a solution until
the first time it hits $\Gamma$, and next setting $X^\varepsilon$
constant after this time. Because of the bad regularity properties of
$\tilde{b}$ and $\sigma$, uniqueness is a difficult problem. We are
only able to prove pathwise uniqueness under the assumption that
$X^{\varepsilon}$ a.s.\ never hits $\Gamma$, and we give in
Sections~\ref{sec:dim1} and~\ref{sec:cond-tau-infinite} explicit
conditions ensuring this assumption and other conditions ensuring the
converse. In section~\ref{sec:LDP}, we prove the main result of this
paper: a large deviation principle for $X^{\varepsilon}$ as
$\varepsilon\rightarrow 0$. Finally, in Section~\ref{sec:exitdomain},
we apply this result to the problem of diffusion exit from an
attracting domain. We obtain a lower bound for the time of exit and we
prove that the exit occurs with high probability near points of the
boundary minimizing the quasi-potential.

\section{Description of the model}
\label{sec:a-b-btilde}

We assume for simplicity that the space of phenotypic traits is
$\RR^d$ for some $d\geq 1$ (this may appear as a restrictive
assumption, however see Remark~\ref{rem:X-cpct} below).  The
coefficients $b$, $\tilde{b}$ and $\sigma\sigma^*=a$ of the
SDE~(\ref{eq:SDE-LDP}) are functions on $\RR^d$, explicitly given in
terms of two biological parameters: the fitness function, and the
mutation law. In this section, we first define these parameters, and
then study their regularity.

\subsection{The fitness function}
\label{sec:fitness}

The function $g(y,x)$ from $\RR^d\times\RR^d$ to $\mathbb{R}$ is the
\emph{fitness function}, which measures the selective advantage (or
disadvantage) of a single mutant individual with trait $y$ in a
population with dominant trait $x$
(see~\cite{metz-geritz-al-96,champagnat-06}). If $g(y,x)>0$ (resp.\
$g(y,x)<0$), then the mutant trait $y$ is selectively advantaged
(resp.\ disadvantaged) in a population of trait $x$. With this in
mind, the fact that the fitness function satisfies
\begin{equation}
  \label{eq:fitness-property}
  g(x,x)=0,\quad \forall x\in\RR^d
\end{equation}
is natural (a mutant trait with trait $x$ is neither advantaged nor
disadvantaged in a population with the same trait).

When $g$ is sufficiently regular, we will denote by $\nabla_1g$ the
gradient of $g(y,x)$ with respect to the first variable $y$, and by
$H_{i,j}g$ the Hessian matrix of $g(y,x)$ with respect to the $i$-th
and $j$-th variables ($1\leq i,j\leq 2$).

We introduce the sets
\begin{align}
  \Gamma & =\{x\in\RR^d:\nabla_1g(x,x)=0\},
  \label{eq:Gamma} \\
  \mbox{and}\quad\forall\alpha>0,\quad
  \Gamma_{\alpha} & =\{x\in\RR^d:d(x,\Gamma)\geq\alpha\mbox{\ and\
  }|x|\leq 1/\alpha\}. \label{eq:Gamma-alpha} 
\end{align}
The points of $\Gamma$ are called \emph{evolutionary singularities}.

We assume that
\begin{description}
\item[\textmd{(H1)}] $g(y,x)$ is ${\cal C}^2$ on $\RR^{2d}$ with
  respect to the first variable $y$, and $\nabla_1g$ and $H_{1,1}g$
  are bounded and Lipschitz on $\RR^{2d}$.
\end{description}


\begin{rem}
  \label{rem:X-cpct} In most biological applications, the trait space
  is a compact subset ${\cal X}$ of $\RR^d$. However, the boundary of
  the trait space usually corresponds to deleterious traits. In other
  words, $g(y,x)\leq 0$ for all $y$ in the boundary of ${\cal
    X}$. Therefore, assuming that the trait space is unbounded is not
  restrictive, since one can extend the fitness function to $\RR^d$ in
  such a way that $g(y,x)\leq 0$ for all $y\not\in{\cal X}$ and
  $x\in\RR^d$. This amounts to add fictive traits, such that
  individuals holding these traits cannot live.
\end{rem}

\subsection{The mutation law}
\label{sec:mutation-law}

The second biological parameter, $p(x,h)dh$, is the law of $h=y-x$,
where $y$ is a mutant trait born from an individual with trait
$x$.
For all $x\in\RR^d$, we assume that this law is absolutely continuous
with respect to Lebesgue's measure and that it is symmetrical with
respect to $0$ for simplicity. This is a very frequent assumption in
adaptive dynamics models (see
e.g.~\cite{dieckmann-law-96,dieckmann-doebeli-99,kisdi-99}).


We also assume that
\begin{description}
\item[\textmd{(H2)}] $p(x,h)dh$ has finite and bounded third-order moment, and
  there exists a measurable function
  $m:\mathbb{R}_+\rightarrow\mathbb{R}_+$ such that
  \begin{equation*}
    \int(\|h\|^2\vee\|h\|^3)m(\|h\|)dh<+\infty, \mbox{\ or
      equivalently\ }\int_{\RR_+}(r^{d+1}\vee r^{d+2})m(r)dr<+\infty,
  \end{equation*}
  where $\|\cdot\|$ is the standard Euclidean norm in $\RR^d$, and for
  any $x,y\in\RR^d$ and $h\in\mathbb{R}^d$,
  \begin{equation}
    \label{eq:H2-1}
    |p(x,h)-p(y,h)|\leq\|x-y\|m(\|h\|) \quad\mbox{and}\quad p(x,h)\leq m(\|h\|).
  \end{equation}
\end{description}
We will denote by (H) the two assumptions (H1) and (H2).

Assumption~(\ref{eq:H2-1}) is satisfied for classical jump measures
taken in applications. For example, it holds when $p(x,h)dh$ is
Gaussian for all $x\in\RR^d$, with covariance matrix $K(x)$ uniformly
non-degenerate, bounded and Lipschitz on $\mathbb{R}^d$.

Assumption~(H2) trivially implies the following property.
\begin{lemma}
  \label{lem:H2}
  Assume~(H2). Let $S=\mathbb{R}^d$ or $S=\{h:h\cdot u>0\}$ for some
  $u\in\mathbb{R}^d\setminus\{0\}$, and let $f$ be a function from
  $\mathbb{R}^d$ to $\mathbb{R}$ such that $f(0)=0$ and
  \begin{equation}
    \label{eq:hyp-f-H2}
    \forall x,y\in\mathbb{R}^d,\ |f(x)-f(y)|\leq
    K\|x-y\|\max\{\|x\|,\|y\|,\|x\|^2,\|y\|^2\}    
  \end{equation}
  for some constant $K$. Then, the function
  $\phi(x)=\int_Sf(h)p(x,dh)$ is globally Lipschitz on $\RR^d$.
\end{lemma}
Note that, in the previous statement, since $f(0)=0$, $|f(h)|\leq
K(\|h\|^2\vee\|h\|^3)$. Thus, the function $\phi$ is well-defined.

As a consequence of this result, (H2) also implies the following
property, needed in the sequel to control the non-degeneracy of the
matrix $a(x)$:
\begin{equation}
  \label{eq:H4}
  \forall\alpha>0,\ \displaystyle{\inf_{\|x\|\leq 1/\alpha,\ u,v\in\mathbb{R}^d,\
      \|u\|=\|v\|=1}\int_{\mathbb{R}^d}|h\cdot u|^2|h\cdot
    v|p(x,h)dh>0},
\end{equation}
where $u\cdot v$ denotes the standard Euclidean inner product between
$u$ and $v\in\RR^d$.  Indeed, $\int_{\mathbb{R}^d}|h\cdot u|^2|h\cdot
v|p(x,h)dh$ is a continuous and positive function of
$(x,u,v)$. Therefore, its minimum on a compact set is positive.

\begin{rem}
  \label{rem:H2}
  Lemma~\ref{lem:H2} is the only consequence of~(H2) that will be used
  below. Assumption~(H2) could be replaced by any condition ensuring
  this result. In particular, it would be sufficient to assume
  regularity of the probability measure $p(x,h)dh$ with respect to
  appropriate Kantorovich
  metrics~\cite{rachev-91}. See~\cite{champagnat-04} for such conditions.
\end{rem}

\subsection{The diffusion model of adaptive dynamics}
\label{sec:diffusion-LDP}

The diffusion model of~\cite{champagnat-ferriere-al-06} is given in
dimension 1. However, the computation of its parameters can be easily
generalized to a multidimensional setting (see~\cite{champagnat-04}
for details). The parameters $b=(b_1,\ldots,b_d)$,
$\tilde{b}=(\tilde{b}_1,\ldots,\tilde{b}_d)$ and
$a=\sigma\sigma^*=(a_{kl})_{1\leq k,l\leq d}$, where $^*$ denotes the
matrix transpose operator, are given by the following expressions: for
all $x\in\RR^d$,
\begin{gather}
  b_k(x)=\int_{\mathbb{R}^d}h_k[\nabla_1g(x,x)\cdot h]_+p(x,h)dh,
  \notag \\
  \tilde{b}_k(x)=\frac{1}{2}\int_{\{ h\cdot\nabla_1g(x,x)>0\}
        }h_k(h^*H_{1,1}g(x,x)h)p(x,h)dh \notag \\
  \mbox{and}\quad a_{kl}(x)=
  \int_{\mathbb{R}^d}h_kh_l[h\cdot\nabla_1g(x,x)]_+p(x,h)dh.
  \label{eq:coefficients-diffusion-LDP}
\end{gather}
We also define
\begin{equation*}
  b^{\varepsilon}=b+\varepsilon\tilde{b}.
\end{equation*}
and the matrix $\sigma$ appearing in~(\ref{eq:SDE-LDP}) as the unique
real symmetrical positive $d\times d$ square root of $a$.

Observe that, for all $x\in\Gamma$, $a(x)=b(x)=\tilde{b}(x)=0$. Thus,
points of $\Gamma$ are possible rest points of solutions
of~(\ref{eq:SDE-LDP}).

The regularity of these parameters is given in the following result.
\begin{prop}
  \label{prop:a-b-btilde}
   Assume (H).
   \begin{description}
   \item[\textmd{(i)}] $a$ and $b$ are globally Lipschitz and bounded
     on $\RR^d$, and $\tilde{b}$ is bounded on $\RR^d$ and locally
     Lipschitz on $\RR^d\setminus\Gamma$.
   \item[\textmd{(ii)}] The matrix $a$ is symmetrical and non-negative
     on ${\cal X}$, $a(x)=0$ if $x\in\Gamma$, and $a(x)$ is positive
     definite if $x\in\RR^d\setminus\Gamma$. For all $\alpha>0$, there
     exists $c>0$ such that $\Gamma_{\alpha}\subset\{
     x\in\RR^d,\:\forall s\in\mathbb{R}^d,\ s^*a(x)s\geq c\|s\|^2\}$,
     where $\Gamma_{\alpha}$ is defined in~(\ref{eq:Gamma-alpha}).
   \item[\textmd{(iii)}] The symmetrical square root $\sigma$ of $a$
     is bounded, H\"older with exponent $1/2$ on $\RR^d$ and locally
     Lipschitz on $\RR^d\setminus\Gamma$.
  \end{description}
\end{prop}

\paragraph{Proof}
In all this proof, the constant $C$ may change from line to line. 

Let us start with Point~(i). The functions $a$, $b$ and $\tilde{b}$
are trivially bounded. Fix $x$ and $y$ in $\RR^d$. For $1\leq k\leq
d$,
\begin{equation*}
  \begin{aligned}
    |b_k(x)-b_k(y)| & \leq\left|\int_{\mathbb{R}^d}h_k
      ([\nabla_1g(x,x)\cdot h]_+-[\nabla_1g(y,y)\cdot
      h]_+)p(x,h)dh\right| \\
    & +\left|\int_{\mathbb{R}^d}h_k[\nabla_1g(x,x)\cdot
      h]_+(p(x,h)-p(y,h))dh\right|.
  \end{aligned}
\end{equation*}
Since $|[a]_+-[b]_+|\leq|a-b|$ and $\nabla_1g$ is Lipschitz, the
first term of the right-hand side is less than $C\|x-y\|M_2$, where
$M_2$ is a bound for the second-order moments of $p(x,h)dh$. Since the
second term is equal to
$$
\left|\int_{\{h\cdot\nabla_1 g(x,x)>0\}}h_k\nabla_1g(x,x)\cdot
  h(p(x,h)-p(y,h))dh\right|,
$$
Lemma~\ref{lem:H2} can be applied to bound this term by
$C\|\nabla_1g(x,x)\|\|x-y\|$.
Since $\nabla_1g$ is bounded, it follows that $b$ is Lipschitz on
$\RR^d$. Similarly, $a$ is Lipschitz on $\RR^d$.

Take $x$ and $y$ in $\RR^d\setminus\Gamma$ and let $S=\{
h\in\mathbb{R}^d:h\cdot\nabla_1g(x,x)>0\}$ and $S'=\{
h:h\cdot\nabla_1g(y,y)>0\}$. We also denote by $S^c$ (resp.\ $S'^c$)
the complement of $S$ (resp.\ $S'$) in $\RR^d$. Then,
\begin{equation}
  \label{eq:continuity-btilde}
  \begin{aligned}
    2|\tilde{b}_k(x)-\tilde{b}_k(y)| &
    \leq\left|\int_{S\cap S'}h_k
      [h^*(H_{1,1}g(x,x)-H_{1,1}g(y,y))h]p(y,h)dh\right| \\
    & \quad +\left|\int_Sh_k(h^*H_{1,1}g(x,x)h)
    (p(x,h)-p(y,h))dh\right| \\
    & \quad +\left|\int_{S\cap S^{\prime c}}
    h_k(h^*H_{1,1}g(x,x)h)p(y,h)dh\right| \\
    & \quad +\left|\int_{S^c\cap
        S'}h_k(h^*H_{1,1}g(y,y)h)p(y,h)dh\right|.
  \end{aligned}
\end{equation}
By Lemma~\ref{lem:H2}, the first two terms of the right-hand side are
both bounded by $C\|x-y\|$ for some constant $C$. The third term
can be bounded by
$$
C\int_{S\cap S'^c}\|h\|^3m(\|h\|)dh.
$$
Making an appropriate spherical coordinates change of variables, this
quantity can be bounded by
$$
C\theta\int_{\RR_+}r^{d+2}m(r)dr\leq C'\theta,
$$
where $\theta$ is the angle between the vectors $\nabla_1g(x,x)$ and
$\nabla_1g(y,y)$.

Now, fix $\alpha>0$. For all $z\in\Gamma_\alpha$, $\nabla_1
g(z,z)\not=0$. Therefore, $\beta:=\inf_{z\in\Gamma_\alpha}\|\nabla_1
g(z,z)\|>0$. Let $K$ be such that $\nabla_1 g(x,x)$ is $K$-Lipschitz
and let $u=\nabla_1 g(x,x)/\|\nabla_1 g(x,x)\|$ and $v=\nabla_1
g(y,y)/\|\nabla_1 g(y,y)\|$. Then
$$
\|u-v\|\leq \frac{\|\nabla_1g(x,x)-\nabla_1 g(y,y)\|}{\|\nabla_1
  g(x,x)\|}+\|\nabla_1 g(y,y)\|\left|\frac{1}{\|\nabla_1
    g(x,x)\|}-\frac{1}{\|\nabla_1 g(y,y)\|}\right|\leq
\frac{K\|x-y\|}{\beta}.
$$
Now, on the one hand $\sin (\theta/2)=\|u-v\|/2$ and on the other
hand, $\sin x\geq (2\sqrt{2}/\pi) x$ for all $0\leq x\leq
\pi/4$. Therefore, for any $x,y\in\Gamma_\alpha$ such that
$\|x-y\|\leq \sqrt{2}\beta/K$, we have
$$
\theta\leq\frac{K\pi}{2\beta}\|x-y\|.
$$
Therefore, for any $x,y\in\Gamma_\alpha$ such that $\|x-y\|\leq
\sqrt{2}\beta/K$,
$$
\left|\int_{S\cap S^{\prime c}}
    h_k(h^*H_{1,1}g(x,x)h)p(y,h)dh\right|\leq C_\alpha\|x-y\|,
$$
where the constant $C_\alpha$ depends only on $\alpha$. Proceeding as
before for the last term of~(\ref{eq:continuity-btilde}), we obtain
that $\tilde{b}$ is uniformly Lipschitz on any convex compact subset
of $\RR^d\setminus\Gamma$, ending the proof of Point~(i).  \medskip

Concerning Point~(ii), $a$ is obviously symmetrical, and for any
$s=(s_1,\ldots,s_d)\in\mathbb{R}^d$, using the symmetry of $p(x,h)dh$,
\begin{equation*}
  \begin{aligned}
    s^*a(x)s & =\int_{\mathbb{R}^d}(h\cdot
    s)^2[h\cdot\nabla_1g(x,x)]_+p(x,h)dh \\
    & =\frac{1}{2}\int_{\mathbb{R}^d}(h\cdot
    s)^2|h\cdot\nabla_1g(x,x)|p(x,h)dh.
  \end{aligned}
\end{equation*}
This is non-negative for all $s\in\mathbb{R}^d$, and is non-zero if
$s\neq 0$ and $x\not\in\Gamma$.

Fix $\alpha>0$, $x\in\Gamma_{\alpha}$, and
$s=(s_1,\ldots,s_d)\in\mathbb{R}^d$. We denote by $u$ and $v$ the unit
vectors of $\mathbb{R}^d$ such that $s=\| s\|u$ and
$\nabla_1g(x,x)=\|\nabla_1g(x,x)\|v$. Then
\begin{equation}
  \label{eq:a-unif-non-degenerate}
  \begin{aligned}
    s^*a(x)s & =\frac{1}{2}\|s\|^2\|\nabla_1g(x,x)\|
    \int_{\mathbb{R}^d}|h\cdot u|^2|h\cdot v|p(x,dh) \\
    & \geq C_{\alpha}\| s\|^2\|\nabla_1g(x,x)\|
  \end{aligned}
\end{equation}
where $C_{\alpha}>0$ by~(\ref{eq:H4}). Since $\Gamma_\alpha$ is a
compact subset of $\RR^d$, we also have
$\inf_{x\in\Gamma_\alpha}\|\nabla_1 g(x,x)\|>0$, completing the proof
of Point~(ii).
\medskip

Finally, Point~(iii) follows from the facts that $a$ is globally
Lipschitz on $\RR^d$ and that the symmetric square root function on
the set of symmetric positive $d\times d$ matrices is globally
$1/2$-H\"older, and Lipschitz in $\{a\in{\cal S}_+:\forall
s\in\mathbb{R}^d,\ s^*as\geq c\|s\|^2\}$ for any $c>0$. A proof of
these facts can be found for example in~\cite{serre-02}.\hfill$\Box$

\section{Strong existence, pathwise uniqueness and strong Markov
  property}
\label{sec:generalities}

Our goal in this section is to construct a particular, strong Markov
solution of the SDE~(\ref{eq:SDE-LDP}), identify the difficulty for
pathwise uniqueness and give some conditions solving this difficulty,
both in the one-dimensional case and the general case.

We fix $\varepsilon>0$ until the end of this section.

\subsection{Strong existence and pathwise uniqueness: construction of
  a particular solution of~(\ref{eq:SDE-LDP})}
\label{sec:construction-X}

\begin{prop}
  \label{prop:exist-LDP}
  Assume~(H). For any filtered probability space $(\Omega,{\cal
    F},({\cal F}_t)_{t\geq 0},\PP,W)$ equipped with a $d$-dimensional
  standard Brownian motion $W$, and for any $x\in\RR^d$, there exists
  a ${\cal F}_t$-adapted process $X^{\varepsilon,x}$ on $\Omega$ a.s.\
  solution of~(\ref{eq:SDE-LDP}) with initial state $x$, such that
  $X^{\varepsilon,x}_t$ is constant after $\tau$, where
  \begin{equation}
    \label{eq:def-tau}
    \tau=\inf\{ t\geq 0:X^{\varepsilon,x}_t\in\Gamma\}.
  \end{equation}
  Moreover, this process is the unique solution of~(\ref{eq:SDE-LDP})
  up to indistinguishability satisfying
  $X^{\varepsilon,x}_t=X^{\varepsilon,x}_\tau$ for all $t\geq\tau$
  a.s.
\end{prop}

\paragraph{Proof}
By Proposition~\ref{prop:a-b-btilde}, the functions $\tilde{b}$ and
$\sigma$ are bounded and locally Lipschitz on
$\RR^d\setminus\Gamma$. Moreover, $b$ is bounded and globally
Lipschitz on $\RR^d$.

Assume that $x\not\in\Gamma$ and fix $\alpha>0$ such that
$x\in\Gamma_\alpha$. Since $\Gamma_\alpha$ is a compact subset of
$\RR^d\setminus\Gamma$. one can construct $\tilde{b}^\alpha$ (resp.\
$\sigma^\alpha$) an extension to $\RR^d$ of $\tilde{b}$ (resp.\
$\sigma$) restricted to $\Gamma_\alpha$ such that $\tilde{b}^\alpha$
(resp.\ $\sigma^\alpha$) is bounded and globally Lipschitz on $\RR^d$
(resp.\ bounded, globally Lipschitz and uniformly non-degenerate on
$\RR^d$). Then, strong existence and pathwise uniqueness for the SDE
\begin{equation}
  \label{eq:SDE-Markov}
  d\tilde{X}_t^{\varepsilon,\alpha}=(b(\tilde{X}_t^{\varepsilon,\alpha})
  +\varepsilon\tilde{b}^{\alpha}(\tilde{X}_t^{\varepsilon,\alpha}))dt
  +\sqrt{\varepsilon}\sigma^{\alpha}(\tilde{X}_t^{\varepsilon,\alpha})dW_t
\end{equation}
with initial condition $\tilde{X}^{\varepsilon,\alpha}_0=x$ are
well-known results. Let
$$
\tau_{\alpha}=\inf\{t\geq 0:\tilde{X}^{\varepsilon,\alpha}_t\not\in\Gamma_{\alpha}\}.
$$
By pathwise uniqueness, for any $\alpha,\alpha'>0$,
$\tilde{X}^{\varepsilon,\alpha}_t=\tilde{X}^{\varepsilon,\alpha'}_t$
for all $t\leq\tau_\alpha\wedge\tau_{\alpha'}$ a.s. Therefore, the
process $X^{\varepsilon,x}$ defined by
$X^{\varepsilon,x}_t=\tilde{X}^{\varepsilon,\alpha}_t$ for
$t\leq\tau_\alpha$ is a solution of~(\ref{eq:SDE-LDP}) for
$t<\sup_{\alpha>0}\tau_\alpha=\tau$.

On the event $\{\tau=+\infty\}$, this gives a strong solution
of~(\ref{eq:SDE-LDP}). On the event $\{\tau<\infty\}$, as a solution
to~(\ref{eq:SDE-LDP}), the semimartingale
$(X^{\varepsilon,x}_t,t<\tau)$ has a uniformly Lipschitz finite
variation part (since $b^\varepsilon$ is bounded), and a local
martingale part which is uniformly in $L^2$, and thus uniformly
integrable, on finite time intervals (since $\sigma$ is
bounded). Therefore, on the event $\{\tau<\infty\}$, the random
variable 
$$
X^{\varepsilon,x}_\tau:=\lim_{t\uparrow
  \tau}X^{\varepsilon,x}_t
$$
is a.s.\ well-defined and finite. Since
$b(x)=\tilde{b}(x)=\sigma(x)=0$ for all $x\in\Gamma$, defining
$X^{\varepsilon,x}_t=X^{\varepsilon,x}_\tau$ for $t\geq\tau$. provides a
strong solution of~(\ref{eq:SDE-LDP}).

In the case where $x\in\Gamma$, setting $X^{\varepsilon,x}_t=x$ for
all $t\geq 0$ trivially provides a strong solution
of~(\ref{eq:SDE-LDP}).

Now, by pathwise uniqueness for~(\ref{eq:SDE-Markov}), there is
pathwise uniqueness for~(\ref{eq:SDE-LDP}) until time
$\tau$. Therefore, the process $X^{\varepsilon,x}$ we constructed
above is the unique solution of~(\ref{eq:SDE-LDP}) constant after time
$\tau$.\hfill$\Box$ \bigskip

The following result is a trivial consequence of the previous one.
\begin{prop}
  \label{prop:uniqueness}
  With the same assumption and notation as in
  Proposition~\ref{prop:exist-LDP}, assume that, for some
  $x\in\RR^d\setminus\Gamma$,
  \begin{equation}
    \label{eq:hyp-tau-finite}
    \PP(X^{\varepsilon,x}_t\not\in\Gamma,\ \forall t\geq 0)=\PP_x(\tau=\infty)=1,
  \end{equation}
  where $\PP_x$ is the law of $X^{\varepsilon,x}$. Then, pathwise
  uniqueness holds for~(\ref{eq:SDE-LDP}) with initial state $x$.
\end{prop}

The question whether pathwise uniqueness also holds for the whole
trajectory when it can hit $\Gamma$ in finite time is
difficult. Because of the singularities of our diffusion ($\tilde{b}$
discontinuous and $\sigma$ degenerate and non-Lipschitz), no standard
technique apply in dimension two or more. In the one-dimensional case,
general criterions of Engelbert and Schmidt exist on pathwise
uniqueness (see~\cite{karatzas-shreve-88}). However, the nature of our
singularity corresponds precisely to a situation where the criterion
does not allow to conclude.
The combination of our singularities is also incompatible with
classical results about uniqueness in law.

Therefore, it is desirable to have conditions
ensuring~(\ref{eq:hyp-tau-finite}) or its converse. This is done is
Sections~\ref{sec:dim1} and~\ref{sec:cond-tau-infinite}. These results
will also be useful in Section~\ref{sec:exitdomain}.

\subsection{Strong Markov property}
\label{sec:uniqueness-Markov}

The strong Markov property for solutions of SDEs is known to be linked
to the uniqueness of solutions to the corresponding martingale
problem. Here, we cannot prove uniqueness in general, but the strong
Markov property can be easily proved.

\begin{prop}
  \label{prop:str-M-pty}
  Assume~(H). Then the family $(X^{\varepsilon,x})_{x\in\RR^d}$ of
  solutions of~(\ref{eq:SDE-LDP}) constructed in
  Proposition~\ref{prop:exist-LDP} satisfy the strong Markov property.
\end{prop}

\paragraph{Proof}
Let $x$ be a fixed point of $\RR^d$, $S$ be a ${\cal F}_t$-stopping
time and $\varphi$ be a bounded and continuous function from $\RR^d$
to $\RR$. We want to prove that
\begin{equation*}
  \EE(\varphi(X^{\varepsilon,x}_{S+t})\mid{\cal F}_S)
  =\EE(\varphi(X^{\varepsilon,x}_{S+t})\mid X^{\varepsilon,x}_S).
\end{equation*}
Since $X^{\varepsilon,x}_t$ is constant after time
$\tau$, this is equivalent to the existence of a
Lebesgue-measurable function $f:\mathbb{R}^d\rightarrow\mathbb{R}$
such that
\begin{equation*}
  \EE(\varphi(X^{\varepsilon,x}_{(S+t)\wedge\tau})\mid {\cal
  F}_S)=f(X^{\varepsilon,x}_S).
\end{equation*}

Recall the definition of $\tau_\alpha$ and
$\tilde{X}^{\varepsilon,\alpha}$ with initial condition $x$ in the
proof of Proposition~\ref{prop:exist-LDP}. Since there is strong
existence and pathwise uniqueness for~(\ref{eq:SDE-Markov}), the
strong Markov property holds for
$\tilde{X}^{\varepsilon,\alpha}$~\cite[Thm.\:5.4.20]{karatzas-shreve-88}. Therefore,
for any $\alpha>0$, there is a bounded Lebesgue-measurable function
$f_{\alpha}$ such that
\begin{equation*}
  \EE(\mathbf{1}_{\tau_{\alpha}>S}
  \varphi(\tilde{X}^{\varepsilon,\alpha}_{(S+t)\wedge\tau_{\alpha}})\mid{\cal F}_S)
  =\mathbf{1}_{\tau_{\alpha}>S}f_{\alpha}(\tilde{X}^{\varepsilon,\alpha}_S).
\end{equation*}
Since $\tilde{X}^{\alpha,\varepsilon}_t=X^{\varepsilon,x}_t$ for all
$t\leq\tau_\alpha$, this yields
\begin{equation*}
  \EE(\mathbf{1}_{\tau_{\alpha}>S}
  \varphi({X}^{\varepsilon,x}_{(S+t)\wedge\tau_{\alpha}})\mid{\cal F}_S)
  =\mathbf{1}_{\tau_{\alpha}>S}f_{\alpha}({X}^{\varepsilon,x}_S).
\end{equation*}
Observing that
$\mathbf{1}_{\tau>S}=\mathbf{1}_{X^{\varepsilon,x}_S\not\in\Gamma}$ is
$\sigma(X^{\varepsilon,x}_S)$-measurable, we deduce that
\begin{equation*}
  \EE(\mathbf{1}_{\tau_{\alpha}>S}
  \varphi({X}^{\varepsilon,x}_{(S+t)\wedge\tau_{\alpha}})\mid{\cal F}_S)+
  \mathbf{1}_{\tau>S\geq\tau_\alpha}f_{\alpha}({X}^{\varepsilon,x}_S)
\end{equation*}
is $\sigma(X^{\varepsilon,x}_S)$-measurable for all
$\alpha>0$. Letting $\alpha$ go to 0, it follows from Lebesgue's
theorem for conditional expectations that this random variable (in
short, r.v.)  a.s.\ converges to $\EE(\mathbf{1}_{\tau>S}
\varphi({X}^{\varepsilon,x}_{(S+t)\wedge\tau})|{\cal F}_S)$. As an
a.s.\ limit of $\sigma(X^{\varepsilon,x}_S)$-measurable r.v., this
r.v.\ is also $\sigma(X^{\varepsilon,x}_S)$-measurable.

Now,
\begin{equation*}
  \EE(\mathbf{1}_{\tau\leq
    S}\varphi(X^{\varepsilon,x}_{(S+t)\wedge\tau})\mid{\cal F}_S)
  =\EE(\mathbf{1}_{\tau\leq S}\varphi(X^{\varepsilon,x}_S)\mid{\cal
    F}_S)
  =\mathbf{1}_{X^{\varepsilon,x}_S\in\Gamma}\varphi(X^{\varepsilon,x}_S),
\end{equation*}
which is also $\sigma(X^{\varepsilon,x}_S)$-measurable. This ends the
proof of Proposition~\ref{prop:str-M-pty}.\hfill$\Box$

\subsection{Study of $\PP(\tau=\infty)$: the case of dimension 1}
\label{sec:dim1}


As we saw above, the uniqueness of $X^{\varepsilon,x}$ relies on the
fact that $\mathbb{P}_x(\tau=\infty)=1$, where $\tau$ has been defined
in~(\ref{eq:def-tau}) and where $\PP_x$ is the law of
$X^{\varepsilon,x}$. Our goal in this section and the following one is to give
conditions under which this is true (or false). 

In this section, we assume that $d=1$. In this case, an elementary
calculation gives the following formulas for $a$, $b$ and $\tilde{b}$:
\begin{equation*}
  \begin{gathered}
    b(x)=\frac{M_2(x)}{2}\partial_1g(x,x),\quad
    \tilde{b}(x)=\frac{M_3(x)}{4}\mbox{sign}[\partial_1g(x,x)]
    \partial^2_{1,1}g(x,x), \\
    a(x)=\frac{M_3(x)}{2}|\partial_1g(x,x)|,\quad\mbox{where}\quad
    M_k(x)=\int_{\mathbb{R}}|h|^kp(x,h)dh \\
    \mbox{and}\quad\mbox{sign}(x)=-1\mbox{\ if\ }x<0;\ 0\mbox{\ if\ }x=0;\
    1\mbox{\ if\ }x>0.
  \end{gathered}
\end{equation*}

In the following result, we use the fact that
$\partial_{1,1}^2g(x,x)+2\partial_{1,2}^2g(x,x)+\partial_{2,2}^2g(x,x)=0$
for all $x\in\RR$, which follows from differentiation
of~(\ref{eq:fitness-property}).

\begin{thm}
  \label{thm:dim1}
  Assume (H). Assume also that $d=1$ and $g$ is ${\cal C}^3$ with
  bounded third-order derivatives. Let $x\not\in\Gamma$ and define
  $c=\sup\{y\in\Gamma,y<x\}$, $c'=\inf\{y\in\Gamma,y>x\}$, and assume
  that $-\infty<c<c'<\infty$,
  $\partial^2_{1,1}g(c,c)+\partial^2_{1,2}g(c,c)\not=0$ and
  $\partial^2_{1,1}g(c',c')+\partial^2_{1,2}g(c',c')\not=0$. We
  can define
  \begin{equation}
    \label{eq:def-alpha-beta}
    \begin{gathered}
      \alpha:=\frac{\partial^2_{1,1}g(c,c)}{\partial^2_{1,1}g(c,c)+
        \partial^2_{1,2}g(c,c)}=\frac{2\partial^2_{1,1}g(c,c)}
      {\partial^2_{1,1}g(c,c)-\partial^2_{2,2}g(c,c)} \\
      \beta:=\frac{\partial^2_{1,1}g(c',c')}{\partial^2_{1,1}g(c',c')+
        \partial^2_{1,2}g(c',c')}=\frac{2\partial^2_{1,1}g(c',c')}
      {\partial^2_{1,1}g(c',c')-\partial^2_{2,2}g(c',c')}.
    \end{gathered}
  \end{equation}
  \begin{description}
  \item[\textmd{(a)}] If $\alpha\geq1$ and $\beta\leq-1$, then
    $\mathbb{P}_x(\tau=\infty)=1$ and the process $X^{\varepsilon,x}$ is
    recurrent in $(c,c')$.
  \item[\textmd{(b)}] If $\alpha\geq1$ and $\beta>-1$, then
    $\mathbb{P}_x(\tau<\infty)=1$ and $\mathbb{P}(\lim_{t\rightarrow
      \tau}X^{\varepsilon,x}_t=c')=1$.
  \item[\textmd{(c)}] If $\alpha<1$ and $\beta\leq-1$, then
    $\mathbb{P}_x(\tau<\infty)=1$ and $\mathbb{P}(\lim_{t\rightarrow
      \tau}X^{\varepsilon,x}_t=c)=1$.
  \item[\textmd{(d)}] If $\alpha<1$ and $\beta>-1$, then
    $\mathbb{P}_x(\tau<\infty)=1$ and \\ $\mathbb{P}(\lim_{t\rightarrow
      \tau}X^{\varepsilon,x}_t=c)=1-\mathbb{P}(\lim_{t\rightarrow
      \tau}X^{\varepsilon}_t=c')\in(0,1)$.
  \end{description}
\end{thm}

\begin{rembr}
  \label{rems:hyp-thm-d=1}
  \begin{itemize}
  \item When $c=-\infty$ or $c'=\infty$, the calculation below depends
    on the precise behaviour of $g$ and $M_k$ near infinity, and no
    simple general result can be stated.
  \item The biological theory of adaptive dynamics gives a
    classification of evolutionary singularities in dimension $d=1$,
    depending on the values of $\partial_{1,1}^2g$ and
    $\partial_{2,2}^2g$ at these points. Here, the condition
    $\alpha\geq 1$ corresponds, when
    $\partial_{1,1}^2g(c,c)-\partial_{2,2}^2g(c,c)>0$, to the case
    $\partial_{1,1}^2g(c,c)+\partial_{2,2}^2g(c,c)\geq 0$, which
    corresponds in the biological terminology (see \emph{e.g.}
    Diekmann~\cite{diekmann-04}) to a converging stable strategy with
    mutual invasibility, which include the evolutionary branching
    condition; and when
    $\partial_{1,1}^2g(c,c)-\partial_{2,2}^2g(c,c)<0$, to the case
    $\partial_{1,1}^2g(c,c)+\partial_{2,2}^2g(c,c)\leq 0$, which
    corresponds biologically to a repelling strategy without mutual
    invasibility.
  \end{itemize}
\end{rembr}

\paragraph{Proof}
We will use the classical method of removal of drift of Engelbert and
Schmidt and the explosion criterion of Feller (see
e.g.~\cite{karatzas-shreve-88}), which can be applied to
$X^{\varepsilon,x}$, considered as a process with value in $(c,c')$
killed when it hits $c$ or $c'$. These methods involve the two
following functions, defined for a fixed $\gamma\in(c,c')$:
\begin{equation}
  \label{eq:def-p-v}
  \begin{gathered}
    p(x)=\int_{\gamma}^x\exp\left[-2\int_{\gamma}^y
      \frac{b^{\varepsilon}(z)dz}{\varepsilon\sigma^2(z)}\right]dy,\
    \forall x\in(c,c'), \\
    \mbox{and}\quad v(x)=\int_{\gamma}^xp'(y)\int_{\gamma}^y\frac{2dz}
    {\varepsilon p'(z)\sigma^2(z)}dy,\ \forall x\in(c,c').
  \end{gathered}
\end{equation}

Then~\cite[pp.\:345--351]{karatzas-shreve-88}, the statements about
the limit of the process $X^{\varepsilon}_t$ when $t\rightarrow\tau$
and about the recurrence of $X^{\varepsilon}$ depend on whether $p(x)$
is finite or not when $x\rightarrow c$ and $c'$, and
the statements about $\tau$ depends on whether $v(x)$ is finite or not
when $x\rightarrow c$ and $c'$.

Let us compute these limits. We will use the standard notation
$f(x)=o(g(x))$ (resp.\ $f(x)=O(g(x))$, resp.\ $f(x)\sim g(x)$) when
$x\rightarrow a$, if $f(x)/g(x)\rightarrow 0$ when $x\rightarrow a$
(resp.\ $|f(x)|\leq C g(x)$ for some constant $C$ in a neighborhood of
$a$, resp.\ $f(x)/g(x)\rightarrow 1$ when $x\rightarrow a$).

\begin{equation}
  \label{eq:dim1-3}
  \frac{b^{\varepsilon}(x)}{\varepsilon\sigma^2(x)}=
  \frac{b^{\varepsilon}(x)}{\varepsilon a(x)}=\frac{M_2(x)}{\varepsilon
    M_3(x)}\mbox{sign}[\partial_1g(x,x)]+
  \frac{1}{2}\frac{\partial^2_{1,1}g(x,x)}{\partial_1g(x,x)},
\end{equation}
so, for $x<y<\gamma$, the quantity inside the exponential appearing in
the definition of $p$ is
\begin{equation*}
  \int_y^{\gamma}\frac{2M_2(z)}{\varepsilon
    M_3(z)}\mbox{sign}[\partial_1g(z,z)]dz+
  \int_y^{\gamma}\frac{\partial^2_{1,1}g(z,z)}{\partial_1g(z,z)}dz.
\end{equation*}

Since $c\not=-\infty$, the first term is bounded for $c<y<\gamma$ (by
Assumption~(H), $M_3$ is positive and continuous on $[c,c']$), so we
only have to study the second term.

When $y\rightarrow c$, an easy calculation gives
\begin{equation*}
  \frac{\partial^2_{1,1}g(z,z)}{\partial_1g(z,z)}=\frac{\alpha}{z-c}+C+o(1),   
\end{equation*}
where $\alpha$ is defined in~(\ref{eq:def-alpha-beta}), and where $C$
is a constant depending on the derivatives of $g$ at $(c,c)$ up to
order $3$. Consequently, when $y\rightarrow c$,
\begin{align}
  \exp\left[-2\int_{\gamma}^y\frac{b^{\varepsilon}(z)dz}
    {\varepsilon\sigma^2(z)}\right] & =\exp\left[C'+o(1)
    +\int_y^{\gamma}\left(\frac{\alpha}{y-c}+C+o(1)\right)dz\right]
  \notag \\ & e^{C''}(y-c)^{-\alpha}, \label{eq:dim1-6}
\end{align}
as $y\rightarrow c$.

Therefore, if $\alpha<1$, $p(c+)>-\infty$, and if $\alpha\geq 1$,
$p(c+)=-\infty$. The same computation gives the same result when
$x\rightarrow c'$, replacing $\alpha$ by $\beta$.

Now let us compute the limit of $v$ at $c$ and $c'$. Since
$p(c'-)=\infty\Rightarrow v(c'-)=\infty$ and $p(c+)=-\infty\Rightarrow
v(c+)=\infty$~\cite[p.\:348]{karatzas-shreve-88}, we only have to
deal with the cases $\alpha<1$ and $\beta>-1$.

Equation~(\ref{eq:dim1-6}) yields $p'(y)\sim e^C(y-c)^{-\alpha}$, so,
for some constant $C$,
\begin{equation*}
  \frac{2}{\varepsilon p'(z)a(z)}\sim C(z-c)^{\alpha-1},
\end{equation*}
since
\begin{equation*}
  a(z)=\frac{M_3(z)}{2}|\partial_1g(z,z)|\sim
  \frac{M_3(c)}{2}|\partial^2_{1,1}g(c,c)+ \partial^2_{1,2}g(c,c)|(z-c).
\end{equation*}

If $\alpha<0$, when $y\rightarrow c$,
$p'(y)\int_y^{\gamma}\frac{2dz}{\varepsilon p'(z)a(z)}\sim
-Cp'(y)(y-c)^{\alpha}$ is bounded on $(c,\gamma)$, and so
$v(c+)<\infty$. If $\alpha=0$,
$p'(y)\int_y^{\gamma}\frac{2dz}{\varepsilon p'(z)a(z)}\sim
C\log(y-c)$, which has a finite integral on $(c,\gamma)$, so
$v(c+)<\infty$. Finally, if $0<\alpha<1$,
$\int_y^{\gamma}\frac{2dz}{\varepsilon p'(z)a(z)}$ is bounded, so
$v(c+)<\infty$ is equivalent to the convergence of the integral
$\int_c^{\gamma}p'(y)dy$, which holds since
$p'(y)\sim\frac{C}{(y-a)^{\alpha}}$ and $\alpha<1$.\hfill$\Box$

\subsection{Study of $\PP(\tau=\infty)$: the general case}
\label{sec:cond-tau-infinite}

Let us turn now to the case $d\geq 2$. The following result gives
conditions under which $\mathbb{P}_x(\tau=\infty)=1$, based on a
comparison of $d(X^{\varepsilon,x},\Gamma)$ with Bessel processes.

\begin{thm}
  \label{thm:d>1}
  Assume (H). Assume also that $g$ is ${\cal C}^2$ on
  $\mathbb{R}^{2d}$ and that the points of $\Gamma$ are isolated. For
  any $y\in\Gamma$, let ${\cal U}_y$ be a neighborhood of $y$ and
  $a^y>0$ and $a_y>0$ two constants such that $a$ is $a^y$-Lipschitz
  on ${\cal U}_y$ and $s^*a(x)s\geq a_y\|s\|^2\|x-y\|$ for all
  $x\in{\cal U}_y$ and $s\in\mathbb{R}^d$. Define also
  \begin{gather*}
    \tilde{b}_y=\inf_{x\in{\cal U}_y\setminus\{y\}}\frac{x-y}{\|x-y\|}
    \cdot\tilde{b}(x) \\
    \mbox{and}\quad\tilde{b}^y=\sup_{x\in{\cal U}_y\setminus\{y\}}
    \frac{x-y}{\|x-y\|} \cdot\tilde{b}(x).
  \end{gather*}
  \begin{description}
  \item[\textmd{(a)}] If for any $y\in\Gamma$,
    $\frac{\tilde{b}_y+da_y/2}{a^y}\geq 1$, then, for all
    $x\not\in\Gamma$, $\mathbb{P}_x(\tau=\infty)=1$ and
    $\mathbb{P}(\lim_{t\rightarrow+\infty}X^{\varepsilon,x}_t\in\Gamma)=0$.
  \item[\textmd{(b)}] If there exists $y\in\Gamma$ such that
    $\frac{\tilde{b}^y+da^y/2}{a_y}<1$, then, for all
    $x\not\in\Gamma$, $\mathbb{P}(\lim_{t\rightarrow\tau}
    X^{\varepsilon,x}_t=y)>0$.
  \end{description}
\end{thm}

Before proving Theorem~\ref{thm:d>1}, let us give some bounds for the
constants involved in this Theorem. This result makes use of the
notation $B(x,r)$ for the open Euclidean ball of $\RR^d$ centered at
$x$ with radius $r$.

\begin{prop}
  \label{prop:d>1}
  Assume (H). Assume also that $g$ is ${\cal C}^2$ on
  $\mathbb{R}^{2d}$ and that the points of $\Gamma$ are isolated. Fix
  $y\in\Gamma$ and $\alpha>0$ such that
  $B(y,\alpha)\cap\Gamma=\{y\}$. Define
  \begin{equation*}
    C=\displaystyle{\inf_{u,v\in\mathbb{R}^d:\|u\|=\|v\|=1}\int|h\cdot
      u|^2|h\cdot v|p(x,h)dh}.
  \end{equation*}
  $C>0$ by~(\ref{eq:H4}). Let $M_3$ be a bound for the third-order
  moment of $p(x,h)dh$ on $B(y,\alpha)$. Let
  $D=H_{1,1}g(y,y)+H_{1,2}g(y,y)$, and denote by $\lambda^y$ (resp.
  $\lambda_y$) the greatest (resp. the smallest) eigenvalue of
  $D^*D$. Suppose also that $D$ is invertible ($\lambda_y>0$).
  Then, for any $\delta>0$ there exists a neighborhood ${\cal U}_y$ of
  $y$ such that, in the statement of Theorem~\ref{thm:d>1}, we can
  take
  \begin{gather*}
    a^y=M_3\sqrt{\lambda^y}+\delta,\quad a_y=C\sqrt{\lambda_y}-\delta, \\
    \tilde{b}^y<\frac{M_3}{2}\|H_{1,1}g(y,y)\|+\delta\quad\mbox{and}\quad
    \tilde{b}_y>-\frac{M_3}{2}\|H_{1,1}g(y,y)\|-\delta.
  \end{gather*}
\end{prop}

\paragraph{Proof}
It follows from the definition~(\ref{eq:coefficients-diffusion-LDP}) of
$\tilde{b}$ that for $x\not=y$,
\begin{equation}
  \label{eq:d>1-1}
  \frac{x-y}{\|x-y\|}\cdot\tilde{b}(x)=\int_{\{\nabla_1g(x,x)\cdot h>0\}}
  \left(\frac{x-y}{\|x-y\|}\cdot h\right)(h^*H_{1,1}g(x,x)h)p(x,h)dh.
\end{equation}
By assumption, the quantity inside the integral can be bounded by
$\|h\|^3[\|H_{1,1}g(y,y)\|+O(\|x-y\|)]p(x,h)$. Therefore,
\begin{align*}
  \frac{x-y}{\|x-y\|}\cdot\tilde{b}(x) & \leq[\|H_{1,1}g(y,y)\|+O(\|x-y\|)]
  \int_{\{\nabla_1g(x,x)\cdot h>0\}}\|h\|^3p(x,h)dh \\ &
  =\frac{M_3}{2} [\|H_{1,1}g(y,y)\|+O(\|x-y\|)].
\end{align*}
This gives the required bounds for $\tilde{b}^y$ and
$\tilde{b}_y$.

It follows from equation~(\ref{eq:a-unif-non-degenerate}) in the proof
of Proposition~\ref{prop:a-b-btilde}, that, for all $s\in\mathbb{R}^d$
and $x\in\mathbb{R}^d$
\begin{equation*}
  C\|s\|^2\|\nabla_1g(x,x)\|\leq s^*a(x)s\leq M_3\|s\|^2\|\nabla_1g(x,x)\|.
\end{equation*}
Considering an orthonormal basis of $\mathbb{R}^d$ in which $D^*D$ is
diagonal, one has $\lambda_y\|v\|^2\leq \|Dv\|^2=v\cdot
D^*Dv\leq\lambda^y\|v\|^2$ for any $v\in\mathbb{R}^d$.  It remains to
observe that $\nabla_1g(x,x)\sim D(x-y)$ when $x\rightarrow y$ to
obtain the required bounds for $a^y$ and $a_y$.\hfill$\Box$

\paragraph{Proof of Theorem~\ref{thm:d>1}}
Fix $y\in\Gamma$. Let us assume for convenience that $y=0$. By
assumption, to this point of $\Gamma$ is associated a neighborhood
${\cal U}_0$ of $0$ and four constants $a_0>0$, $a^0>0$, $\tilde{b}_0$
and $\tilde{b}^0$. Let $\rho$ be small enough for
$B(\rho):=\{x\in\mathbb{R}^d:\|x\|\leq\rho\}\subset{\cal U}_0$ and
$\Gamma\cap B(\rho)=\{0\}$, and define $\tau_{\rho}:=\inf\{t\geq
0:\|X^{\varepsilon}_t\|=\rho\}$ and $\tau_0=\inf\{t\geq
0:X^{\varepsilon}_t=0\}$, where we omitted the dependence of
$X^{\varepsilon,x}$ with respect to the initial condition. Recall also
the notation $\PP_x$ for the law of $X^\varepsilon$ when
$X^\varepsilon_0=x$.

Theorem~\ref{thm:d>1} can be deduced from the next lemma.
\begin{lembr}
  \label{lem:d-geq-2}
  \begin{description}
  \item[\textmd{(a)}] If $\frac{\tilde{b}_0+da_0/2}{a^0}\geq 1$, then,
    for all $x\in B(\rho)\setminus\{0\}$,
    $\mathbb{P}_x(\tau_{\rho}\leq\tau_0)=1$.
  \item[\textmd{(b)}] If $\frac{\tilde{b}^0+da^0/2}{a_0}<1$, then,
    there exists a constant $c>0$ such that, for all $x\in
    B(\rho/2)\setminus\{0\}$,
    $\mathbb{P}_x(\{\tau_0<\tau_{\rho}\}\cup\{\tau_0=\tau_{\rho}=\infty\mbox{\ 
      and\ }\lim_{t\rightarrow+\infty}X^{\varepsilon}_t=0\})\geq c$.
  \end{description}
\end{lembr}
Together with the strong Markov property of
Proposition~\ref{prop:str-M-pty}, Point~(a) of this lemma easily
implies Theorem~\ref{thm:d>1}~(a), and part~(b) implies
Theorem~\ref{thm:d>1}~(b) if we can prove that for any
$x\in\RR^d\setminus\Gamma$,
$\mathbb{P}_x(\tau_{\rho/2}<\infty)>0$. This can be done as follows.

Let $D$ be any connected bounded open domain $D$ with smooth boundary
containing $B(\rho/2)$. The process $\tilde{X}^{\varepsilon,\alpha}$
of the proof of Proposition~\ref{prop:exist-LDP} has smooth and
uniformly non-degenerate coefficients. Therefore, it is standard to
prove that such a process hits $B(\rho/2)$ before hitting $\partial D$
with positive probability, starting from any $y\in D$. (This may be
proved for example by applying Feynman-Kac's formula to obtain the
elliptic PDE solved by this probability in $D\setminus B(\rho/2)$, and
next applying the strong maximum principle to this PDE.) Choosing
$\alpha$ and $D$ such that $x\in D$ and $D\setminus B(\rho/2)\subset
\Gamma_{\alpha}$, we easily obtain the required estimate.\hfill$\Box$
\bigskip

Before coming to the proof of Lemma~\ref{lem:d-geq-2}, we need to
introduce a few notation: it follows from It\^o's formula that, for
all $t<\tau$,
\begin{multline*}
  \|X^{\varepsilon}_t\| =\|x\|+\int_0^t\frac{1}
  {\|X^{\varepsilon}_s\|}\biggl[X^{\varepsilon}_s\cdot
  (b(X^{\varepsilon}_s)+\varepsilon\tilde{b}(X^{\varepsilon}_s)) \\
  +\frac{\varepsilon}{2}\mbox{Tr}(a(X^{\varepsilon}_s))-
  \frac{\varepsilon}{2}\frac{(X^{\varepsilon}_s)^*}
  {\|X^{\varepsilon}_s\|}a(X^{\varepsilon}_s)
  \frac{X^{\varepsilon}_s}{\|X^{\varepsilon}_s\|}\biggr]ds+M_t,
\end{multline*}
where $\mbox{Tr}$ is the trace operator on $d\times d$ matrices, and
where, for $t<\tau$,
\begin{equation*}
  M_t:=\sqrt{\varepsilon}\int_0^t
  \frac{(X^{\varepsilon}_s)^*}{\|X^{\varepsilon}_s\|}
  \sigma(X^{\varepsilon}_s)dW_s.
\end{equation*}
Let us extend $M_t$ to $t\geq\tau$ by setting $M_t=M_{t\wedge\tau}$
for all $t\geq 0$. Since $\sigma$ is bounded, $M_t$ is a
$\mathbb{L}^2$-martingale in $\mathbb{R}$ with quadratic variation
\begin{equation}
  \label{eq:quadr-M}
  \langle M\rangle_t=\varepsilon\int_0^{t\wedge\tau}
  \frac{(X^{\varepsilon}_s)^*}{\|X^{\varepsilon}_s\|}
  a(X^{\varepsilon}_s)\frac{X^{\varepsilon}_s}{\|X^{\varepsilon}_s\|}ds.
\end{equation}
It follows from Dubins-Schwartz's Theorem (see
e.g.~\cite{karatzas-shreve-88}) that for any $t\geq 0$,
$M_t=B_{\langle M\rangle_t}$, where $B$ is a one-dimensional Brownian
motion.

Define the time change $T_t=\inf\{s\geq 0:\langle M\rangle_s>t\}$ for
all $t\geq 0$. If $t<\langle M\rangle_{\infty}:=\lim_{t\rightarrow
  +\infty}\langle M\rangle_t=\langle M\rangle_\tau$, then $T_t<\infty$
and $\langle M\rangle_{T_t}=t$. For $t<\langle M\rangle_{\infty}$,
define $Y_t=X^{\varepsilon}_{T_t}$.  An easy change of variable shows
that for $t<\langle M\rangle_{\infty}$, $Y_t\not\in\Gamma$, and
\begin{equation*}
  \|Y_t\|=\|x\|+\int_0^tc(Y_s)ds+B_t,
\end{equation*}
where
\begin{equation*}
  c(z)=\|z\|\frac{z\cdot(b(z)+\varepsilon
    \tilde{b}(z))+\varepsilon\mbox{Tr}(a(z))/2}{\varepsilon z^*a(z)z}
  -\frac{1}{2\|z\|}.
\end{equation*}
Using the constants defined in the statement of Theorem~\ref{thm:d>1},
the fact that $b$ is $K$-Lipschitz on $\mathbb{R}^d$, and the fact
that $\mbox{Tr}(a)=\sum_{i=1}^de_i^*ae_i$, where $e_i$ is the
$i^{\mbox{\footnotesize{th}}}$ vector of the canonical basis of
$\mathbb{R}^d$, one easily obtains that, for all $z\in{\cal U}_0$,
\begin{equation*}
  c_1(\|z\|)<c(z)<c_2(\|z\|),
\end{equation*}
where, for $u>0$,
\begin{align*}
  c_1(u) & =\left(\frac{da_0/2+\tilde{b}_0}{a^0}-\frac{1}{2}\right)
  \frac{1}{u}-\frac{2K}{\varepsilon a_0} \\ \mbox{and}\quad
  c_2(u) & =\left(\frac{da^0/2+\tilde{b}^0}{a_0}-\frac{1}{2}\right)
  \frac{1}{u}+\frac{2K}{\varepsilon a_0}.
\end{align*}
Define also the processes $Z^1$ and $Z^2$ strong solutions in
$(0,\infty)$ to the SDEs
\begin{equation*}
  Z^i_t=\|x\|+\int_0^tc_i(Z^i_s)ds+B_t
\end{equation*}
for $i=1,2$, and stopped when they reach $0$. As strong solutions,
these processes can be constructed on the same probability space than
$X^{\varepsilon}$ (and $Y$). Finally, define for $1\leq i\leq 2$ the
stopping times
\begin{gather*}
  \theta^i_0=\inf\{t\geq 0:Z^i=0\} \\
  \mbox{and}\quad\theta^i_{\rho}=\inf\{t\geq 0:Z^i=\rho\}.
\end{gather*}

The proof of Lemma~\ref{lem:d-geq-2} relies on the following three
lemmas. The first one is a comparison result between $Z^1$, $Z^2$ and
$Y$.
\begin{lemma}
  \label{lem:Z1<Y<Z2}
  Almost surely, $Z^1_t\leq\|Y_t\|$ for all  $t<\theta^1_{\rho}\wedge\langle
  M\rangle_{\infty}$, and $\|Y_t\|\leq
  Z^2_t$ for all $t<\theta^2_{\rho}\wedge\langle M\rangle_{\infty}$.
\end{lemma}
The processes $Z^1$ and $Z^2$ are Bessel processes with additional
drifts. The second lemma examines whether these processes hit $0$ in
finite time or not.
\begin{lembr}
  \label{lem:Bessel}
  \begin{description}
  \item[\textmd{(a)}] $Z^1$ is recurrent in $(0,+\infty)$ if and only
    if $\frac{\tilde{b}_0+da_0/2}{a^0}\geq 1$.
  \item[\textmd{(b)}] Let ${\mathbf{P}}_u$ be the law of $Z^2$
    with initial state $u>0$. If $\frac{\tilde{b}^0+da^0/2}{a_0}<1$,
    then, for any $u<\rho$,
    ${\mathbf{P}}_u(\theta^2_0<\theta^2_{\rho})>0$.
  \end{description}
\end{lembr}
The last lemma states that, when $\langle M\rangle_{\infty}<\infty$,
$X^{\varepsilon}$ reaches $\Gamma$ in finite or infinite time.
\begin{lemma}
  \label{lem:lim-Xeps}
  $\{\langle
  M\rangle_{\infty}<\infty\}\subset\{\tau<\infty\}\cup\{\tau=\infty\mbox{\ 
    and\ }\lim_{t\rightarrow+\infty}X^{\varepsilon}_t\in\Gamma\}$ a.s.
\end{lemma}

\paragraph{Proof of Lemma~\ref{lem:d-geq-2}}
Assume first that $\frac{\tilde{b}_0+da_0/2}{a^0}\geq 1$, and fix
$x\in B(\rho)\setminus\{0\}$. Then, by Lemma~\ref{lem:Bessel}~(a),
$\theta^1_0=\infty$ and $\theta^1_{\rho}<\infty$ a.s. Moreover, by
Lemma~\ref{lem:Z1<Y<Z2}, for all $t<T_{\theta^1_{\rho}}$,
$\|X^{\varepsilon}_t\|=\|Y_{\langle M\rangle_t}\|\geq Z^1_{\langle
  M\rangle_t}$.

Then, $\langle M\rangle_{\infty}=\infty$ implies a.s.\ that there
exists $t\geq 0$ such that $\langle M\rangle_t=\theta^1_\rho$ and thus
$\tau_{\rho}<\tau_0$. Conversely, by Lemma~\ref{lem:lim-Xeps},
$\langle M\rangle_{\infty}<\infty$ implies a.s.\ that
$\lim_{t\rightarrow\tau}X^{\varepsilon}_t \in\Gamma\setminus\{0\}$,
and thus that $\tau_{\rho}<\tau_0$. This completes the proof of
Lemma~\ref{lem:d-geq-2}~(a).

Now, assume that $\frac{\tilde{b}^0+da^0/2}{a_0}<1$ and fix $x\in
B(\rho/2)$. By Lemma~\ref{lem:Z1<Y<Z2}, for all
$t<T_{\theta^2_{\rho}}$, $\|X^{\varepsilon}_t\|=\|Y_{\langle
  M\rangle_t}\|\leq Z^2_{\langle M\rangle_t}$.

Then, on the event $\{\theta^2_0<\theta^2_{\rho}\}$, $\langle
M\rangle_{\infty}=\infty$ implies a.s.\ that
$\tau_0<\tau_{\rho}$. Conversely, on the event $\{\theta^2_0<\theta^2_{\rho}\}$,
by Lemma~\ref{lem:lim-Xeps}, $\langle M\rangle_{\infty}<\infty$
implies a.s.\ that $\lim_{t\rightarrow\tau}X^{\varepsilon}_t=0$ (where
$\tau$ may be finite or infinite), and thus that $\tau_0<\tau_{\rho}$
or that $\tau_0=\tau_{\rho}=\infty\mbox{\ and\ 
}\lim_{t\rightarrow+\infty}X^{\varepsilon}_t=0$. Hence,
\begin{equation*}
  \mathbb{P}_x(\{\tau_0<\tau_{\rho}\}\cup
  \{\tau_0=\tau_{\rho}=\infty\mbox{\ and\ }\lim_{t\rightarrow+\infty}
  X^{\varepsilon}_t=0\})\geq{\mathbf{P}}_{\|x\|}
  (\theta^2_0<\theta^2_{\rho}).
\end{equation*}
But, applying the Markov property to $Z^2$,
${\mathbf{P}}_{\|x\|}(\theta^2_0<\theta^2_{\rho})
\geq{\mathbf{P}}_{\rho/2}(\theta^2_0<\theta^2_{\rho})$ for any $x\in
B(\rho/2)$. Since this is positive by Lemma~\ref{lem:Bessel}~(b), the
proof of Lemma~\ref{lem:d-geq-2}~(b) is completed.\hfill$\Box$

\paragraph{Proof of Lemma~\ref{lem:Z1<Y<Z2}}
First, remind that $Y_t$ is defined only for $t<\langle
M\rangle_{\infty}$. Observe that for $t<\theta^1_0\wedge\langle
M\rangle_{\infty}$,
\begin{equation*}
  \|Y_t\|-Z^1_t=\int_0^t(c(Y_s)-c_1(Z^1_s))ds.
\end{equation*}
If there exists $t_0<\theta^1_{\rho}\wedge\theta^1_0\wedge\langle
M\rangle_{\infty}$ such that $\|Y_{t_0}\|=Z^1_{t_0}$, then
$(\|Y\|-Z^1)'(t_0)=c(Y_{t_0})-c_1(Z^1_{t_0})=c(Y_{t_0})-c_1(\|Y_{t_0}\|)>0$,
and therefore, $\|Y_t\|>Z^1_t$ for $t>t_0$ in a neighborhood of $t_0$.
Consequently, $Z^1_t\leq\|Y_t\|$ for any
$t<\theta^1_{\rho}\wedge\theta^1_0\wedge\langle M\rangle_{\infty}$.
Since $Z^1_t=0$ for $t\geq\theta^1_0$, this inequality actually holds
for $t<\theta^1_{\rho}\wedge\langle M\rangle_{\infty}$. The proof of
the other inequality is similar.\hfill$\Box$

\paragraph{Proof of Lemma~\ref{lem:Bessel}}
The proof relies on the same functions $p$ and $v$ as in the proof of
Theorem~\ref{thm:dim1}. They are given by~(\ref{eq:def-p-v}), where
$b^{\varepsilon}$ has to be replaced by $c_i$, and $\varepsilon a$ by
$1$. For the process $Z^1$, if we fix $\gamma>0$, then, for any $x>0$,
\begin{align*}
  p(y) & =\int_{\gamma}^y\exp\left[-2\int_{\gamma}^u c_1(z)dz\right]du
  =-\int^{\gamma}_y\exp\left[2k\int^{\gamma}_u
    \frac{dz}{z}-k'(\gamma-u)\right]du \\ & =-C\int_y^{\gamma}u^{-2k}e^{k'u}du,
\end{align*}
where we have used the constants $k=\frac{\tilde{b}_0+da_0/2}{a^0}
-\frac{1}{2}$ and $k'=\frac{4K}{\varepsilon a_0}$.  Consequently,
$p(0+)=-\infty$ if and only if $2k\geq 1$, and $p(+\infty)=+\infty$,
which yields~(a). A similar computation for $Z^2$ gives that
$p(0+)>-\infty$ if and only if $\frac{\tilde{b}^0+da^0/2}{a_0}<1$,
which yields Lemma~\ref{lem:Bessel}~(b).\hfill$\Box$

\paragraph{Proof of Lemma~\ref{lem:lim-Xeps}} Assume that
$\mathbb{P}(\{\langle M\rangle_{\infty}<\infty\}\cap
\{\lim_{t\rightarrow+\infty}X^{\varepsilon}_t\in\Gamma\}^c)>0$. Then,
there exists $\alpha>0$ such that
\begin{equation*}
  \delta:=\mathbb{P}(\langle M\rangle_{\infty}<\infty,\
  \limsup_{t\rightarrow+\infty}d(X^{\varepsilon}_t,\Gamma)\geq\alpha)>0.
\end{equation*}
Define for any $t>0$ the stopping time $\tau_{\alpha,t}=\inf\{s\geq
t:d(X^{\varepsilon}_s,\Gamma)\geq\alpha\}$. Then , for any $t>0$,
\begin{equation}
  \label{eq:contradiction}
  \mathbb{P}(\langle M\rangle_{\infty}<\infty,\
  \tau_{\alpha,t}<\infty)\geq\delta.
\end{equation}
We will obtain a contradiction from this statement thanks to the
following inequality:
for any $\varepsilon<1$, $h\in(0,1)$ and any stopping time $S$ a.s.\
finite,
\begin{equation*}
  \mathbb{E}\left[\sup_{0<u<h}\|X^{\varepsilon}_{S+u}-
    X^{\varepsilon}_S\|^2\right]\leq 10C^2h,
\end{equation*}
where $C$ is a bound for $b$, $\tilde{b}$ and $\sigma$ on $\RR^d$.
This is a straightforward consequence of the inequality
\begin{equation*}
  \|X^{\varepsilon}_{S+u}-X^{\varepsilon}_S\|^2\leq
  2\left(\int_S^{S+u}\|b(X^{\varepsilon}_s)
    +\varepsilon\tilde{b}(X^{\varepsilon}_s)\|ds\right)^2
  +2\sqrt{\varepsilon}\left\|\int_S^{S+u}
    \sigma(X^{\varepsilon}_s)dW_s\right\|^2
\end{equation*}
and of Doob's inequality.

Taking $h=\delta\alpha^2/80C^2$ and $S=\tau_{\alpha,t}\wedge T$, we
get
\begin{equation*}
  \mathbb{P}\left(\sup_{0<u<h}\|X^{\varepsilon}_{\tau_{\alpha,t}\wedge T+u}-
    X^{\varepsilon}_{\tau_{\alpha,t}\wedge T}\|
    >\frac{\alpha}{2}\right)\leq\frac{\delta}{2}.
\end{equation*}
Letting $T\rightarrow+\infty$,
\begin{equation*}
  \mathbb{P}\left(\tau_{\alpha,t}<\infty,\
    \sup_{0<u<h}\|X^{\varepsilon}_{\tau_{\alpha,t}+u}-
    X^{\varepsilon}_{\tau_{\alpha,t}}\|>\frac{\alpha}{2}\right)
  \leq\frac{\delta}{2}. 
\end{equation*}
Together with inequality~(\ref{eq:contradiction}), this yields the
first line of the following inequality, and the last line makes use
of~(\ref{eq:quadr-M}) and a constant $C>0$ such that $s^*a(x)s\geq
C\|s\|^2$ for any $s\in\mathbb{R}^d$ and $x\in\Gamma_{\alpha/2}$.
\begin{align*}
  \frac{\delta}{2} & \leq\mathbb{P}\left(\langle
    M\rangle_{\infty}<\infty,\ \sup_{0<u<h}\|
    X^{\varepsilon}_{\tau_{\alpha,t}+u}-
    X^{\varepsilon}_{\tau_{\alpha,t}}\|\leq\frac{\alpha}{2}\right) \\
  & \leq\mathbb{P}\left(\langle M\rangle_{\infty}<\infty,\ 
    \inf_{0<u<h}\|X^{\varepsilon}_{\tau_{\alpha,t}+u}\|
    \geq\frac{\alpha}{2}\right) \\
  & \leq\mathbb{P}\left(\langle M\rangle_{\infty}<\infty,\ \langle
    M\rangle_{\tau_{\alpha,t}+h}-\langle
    M\rangle_{\tau_{\alpha,t}}\geq\varepsilon Ch\right).
\end{align*}
Therefore,
\begin{equation*}
  \mathbb{P}\left(\langle M\rangle_{\infty}<\infty,\ \langle
  M\rangle_{\infty}-\langle M\rangle_t\geq\varepsilon
  Ch\right)\geq\frac{\delta}{2}
\end{equation*}
holds for any $t>0$, which is impossible.\hfill$\Box$

\section{Large deviations for $X^{\varepsilon}$ as
  $\varepsilon\rightarrow 0$}
\label{sec:LDP}

Our large deviation result will be obtained by a transfer technique
to carry the LDP from the family
$\{\sqrt{\varepsilon}W\}_{\varepsilon>0}$, where $W$ is a standard
$d$-dimensional Brownian motion (Schilder's Theorem,
e.g.~\cite[p.\:185]{dembo-zeitouni-93}) to the family
$\{X^{\varepsilon}\}_{\varepsilon>0}$, where $X^{\varepsilon}$ is the
solution to the SDE~(\ref{eq:SDE-LDP}) defined in
section~\ref{sec:construction-X}. The method of the proof, adapted
from Azencott~\cite{azencott-80}, consists in constructing a function
$S$ mapping (in some sense) the paths of $\sqrt{\varepsilon}W$ to the
paths of $X^{\varepsilon}$.

\subsection{Statement of the result}
\label{sec:statement-LDP}

We denote by ${\cal C}([0,T],\RR^d)$ (resp.\ ${\cal
  C}^{ac}([0.T],\RR^d)$\:) the set of continuous (resp.\ absolutely
continuous) functions from $[0,T]$ to $\RR^d$.  Fix $T>0$ and
$x\in\RR^d$, and define
\begin{equation*}
  \begin{gathered}
    \forall\psi\in{\cal C}\left([0,T],\RR^d\right),
    \quad t_{\psi}=\inf\{t\in[0,T]:\psi(t)\in\Gamma\}\wedge T \\
    \mbox{and}\quad\tilde{{\cal C}}^{ac}_x([0,T],\RR^d)=
    \{\psi\in{\cal C}^{ac}([0,T],\RR^d)\mbox{\ 
      constant on\ }[t_{\psi},T]\mbox{\ such that\ }\psi(0)=x\}.
  \end{gathered}
\end{equation*}
Then, we define for $\psi\in{\cal C}([0,T],\RR^d)$
\begin{equation}
  \label{eq:I}
  I_{T,x}(\psi)=\left\{
    \begin{array}{l}
      \displaystyle{\frac{1}{2}\int_0^{t_{\psi}}
        [\dot{\psi}(t)-b(\psi(t))]^*a^{-1}(\psi(t))
        [\dot{\psi}(t)-b(\psi(t))]dt} \\ \phantom{+\infty}
      \qquad\qquad\qquad\qquad\qquad\quad\mbox{if}\quad
      \psi\in\tilde{\cal C}^{ac}_x([0,T],\RR^d) \\
      +\infty\qquad\qquad\qquad\qquad\qquad\quad\mbox{otherwise.}
    \end{array}\right.
\end{equation}
By Proposition~\ref{prop:a-b-btilde}~(ii), the inverse matrix
$a^{-1}(x)$ of $a(x)$ is well-defined, symmetric and non-negative for
all $x\not\in\Gamma$, so $I_{T,x}(\psi)$ is well-defined and belongs
to $\mathbb{R}_+\cup\{ +\infty\}$.  When $t_{\psi}=T$, $I_{T,x}(\psi)$
takes the classical form of rate functions for diffusion processes.

This original form of rate function will appear naturally in the
proof. However, as shown in Proposition~\ref{prop:I-not-lsc} below,
this function is \emph{not} lower semicontinuous. Therefore, it is
natural to introduce for all $\psi\in {\cal C}([0,T],\RR^d)$
\begin{equation}
  \label{eq:rate-function}
  \tilde{I}_{T,x}(\psi)=\liminf_{\tilde{\psi}\rightarrow\psi}
  I_{T,x}(\tilde{\psi}),
\end{equation}
which is the biggest lower semicontinuous function on ${\cal
  C}([0,T],\RR^d)$ smaller than $I_{T,x}$.

\begin{thm}
  \label{thm:true-LDP}
  Assume~(H). Assume also that the points of $\Gamma$ are isolated in
  $\mathbb{R}^d$. Fix $T>0$. Then, for any $x\in\RR^d$
  and any open subset $O$ of ${\cal C}([0,T],\RR^d)$,
  \begin{equation}
    \label{eq:true-LB}
    \liminf_{\varepsilon\rightarrow
      0,y\rightarrow x}\varepsilon\ln\mathbb{P}(X^{\varepsilon,y}\in
    O)\geq -\inf_{\psi\in O}\tilde{I}_{T,x}(\psi),
  \end{equation}
  and for any $x\not\in\Gamma$ and any closed subset $C$ of ${\cal
    C}([0,T],\RR^d)$,
  \begin{equation}
    \label{eq:true-UB}
    \limsup_{\varepsilon\rightarrow
      0,y\rightarrow x} \varepsilon\ln\mathbb{P}(X^{\varepsilon,y}\in
    C)\leq -\inf_{\psi\in C}\tilde{I}_{T,x}(\psi).
  \end{equation}
\end{thm}

The general form of the lower and upper bounds~(\ref{eq:true-LB})
and~(\ref{eq:true-UB}) (where the limit is taken over $y\rightarrow
x$) will be useful in Section~\ref{sec:exitdomain}. This general form
requires the restriction that $x\not\in\Gamma$ for the upper bound for
technical reasons. However, this result implies that the following
standard form of LDP holds without any restriction.

\begin{cor}
  \label{cor:true-LDP}
  Assume the conditions of Theorem~\ref{thm:true-LDP}. Then, for any
  $x\in\RR^d$, for any open $O\subset{\cal C}([0,T],\RR^d)$, and for any
  closed $C\subset{\cal C}([0,T],\RR^d)$,
  \begin{gather}
    \label{eq:simple-LB}
    \liminf_{\varepsilon\rightarrow
      0}\varepsilon\ln\mathbb{P}(X^{\varepsilon,x}\in O)\geq 
    -\inf_{\psi\in O}\tilde{I}_{T,x}(\psi), \\
    \label{eq:simple-UB}
    \limsup_{\varepsilon\rightarrow
      0}\varepsilon\ln\mathbb{P}(X^{\varepsilon,x}\in
    C)\leq -\inf_{\psi\in C}\tilde{I}_{T,x}(\psi).
  \end{gather}
\end{cor}

\begin{proof}
  The lower bound~(\ref{eq:simple-LB}) is a trivial consequence
  of~(\ref{eq:true-LB}) and the upper bound~(\ref{eq:simple-UB}) for
  $x\not\in\Gamma$ also trivially follows from~(\ref{eq:true-UB}). If
  $x\in\Gamma$, let us denote by $x$ the constant function of ${\cal
    C}([0,T],\RR^d)$ equal to $x$. In this case,
  $X^{\varepsilon,x}_t=x$ for all $t\geq 0$. Therefore,
  $\PP(X^{\varepsilon,x}\in C)$ equals 1 if the function $x$ belongs
  to $C$, and equals 0 otherwise. Since $\tilde{I}_{T,x}(x)\leq
  I_{T,x}(x)=0$, the upper bound~(\ref{eq:simple-UB}) is clear when
  $x\in\Gamma$.
\end{proof}

\begin{rem}
  \label{rem:canonical-eq}
  As usual for large deviation principles,
  Corollary~\ref{cor:true-LDP} implies the convergence in probability
  of $X^{\varepsilon,x}$ to the solution with initial state $x$ of the
  deterministic ODE
  \begin{equation*}
    \dot{\phi}=b(\phi)
  \end{equation*}
  as $\varepsilon\rightarrow 0$. This ODE is known as the
  \emph{canonical equation of adaptive
    dynamics}~\cite{dieckmann-law-96, champagnat-ferriere-al-01,
    champagnat-meleard-08}.
\end{rem}

In Section~\ref{sec:exitdomain}, we will use the following classical
consequence of Theorem~\ref{thm:LDP}, which can be proved exactly as
Corollary~5.6.15 of~\cite{dembo-zeitouni-93}:
\begin{cor}
  \label{cor:LDP}
  Assume the conditions of Theorem~\ref{thm:true-LDP}. Then, for any
  compact set $K\subset\RR^d$ and for any open $O\subset{\cal
    C}([0,T],\RR^d)$,
  \begin{equation*}
    \liminf_{\varepsilon\rightarrow 0}\varepsilon\ln\inf_{y\in
      K}\mathbb{P}(X^{\varepsilon,y}\in O)\geq-\sup_{y\in
      K}\inf_{\psi\in O}\tilde{I}_{T,y}(\psi),
  \end{equation*}
  and if $K\cap\Gamma=\emptyset$, for any closed $C\subset{\cal
    C}([0,T],\RR^d)$, 
  \begin{equation*}
    \limsup_{\varepsilon\rightarrow
    0}\varepsilon\ln\sup_{y\in K}\mathbb{P}(X^{\varepsilon,y}\in
  C)\leq-\inf_{y\in K,\ \psi\in C}\tilde{I}_{T,y}(\psi).
  \end{equation*}
\end{cor}
\bigskip

We end this subsection with some remarks on the rate functions we
obtain and their links with the classical form of rate functions for
diffusion processes with small noise.

\begin{prop}
  \label{prop:I-not-lsc}
  Assume the conditions of Theorem~\ref{thm:true-LDP}. Assume also
  that there exists an isolated point $y$ of $\Gamma$ such that $g$ is
  ${\cal C}^2$ at $(y,y)$, and that $H_{1,1}g(y,y)+H_{1,2}g(y,y)$ is
  invertible.  Then, for any $x\not\in\Gamma$ and $T>0$, $I_{T,x}$ is
  \emph{not} lower semicontinuous.
\end{prop}
We postpone the proof of this result at the end of this subsection.
\medskip

General large deviation estimates are known for diffusions in $\RR^d$
with small noise using different techniques. For example, Dupuis,
Ellis and Weiss~\cite{dupuis-ellis-al-91} have obtained upper bounds
under very general assumptions. We could have applied their result in
our case (with some modifications since they consider a drift that
does not depend on $\varepsilon$, see Remark~1.2
in~\cite{dupuis-ellis-al-91}) to obtain a large deviations upper bound
with lower semicontinuous rate function
\begin{equation*}
  \hat{I}_{T,x}(\psi)=\left\{
    \begin{array}{l}
      \displaystyle{\frac{1}{2}\int_0^T\mathbf{1}_{\psi(t)\not\in\Gamma}
        [\dot{\psi}(t)-b(\psi(t))]^*a^{-1}(\psi(t))
        [\dot{\psi}(t)-b(\psi(t))]dt} \\ \phantom{+\infty}
      \qquad\qquad\qquad\qquad\qquad\quad\mbox{if\ }
      \psi\in{\cal C}^{ac}([0,T],\RR^d)\mbox{\ and\ }\psi(0)=x \\
      +\infty\qquad\qquad\qquad\qquad\qquad\quad\mbox{otherwise,}
    \end{array}\right.
\end{equation*}
for all $\psi\in{\cal C}([0,T],\RR^d)$. 

Since obviously $\hat{I}_{T,x}\leq I_{T,x}$ and $\hat{I}_{T,x}$ is
lower semicontinuous, we have
$\hat{I}_{T,x}\leq\tilde{I}_{T,x}$. Since $\tilde{I}_{T,x}\leq
I_{T,x}$, this immediately implies that $\hat{I}_{T,x}$ and
$\tilde{I}_{T,x}$ coincide on ${\cal
  C}([0,T],\RR^d\setminus\Gamma)$. Unfortunately, because of the
degeneracy of $a$ on $\Gamma$, we are not able to obtain an explicit
expression for $\tilde{I}_{T,x}(\psi)$ when $\psi_t\in\Gamma$ for some
$t\in[0,T]$. However, it is possible to find simple examples where
these two rate function are not equal: Assume that $d=1$ and $0$ is an
isolated point of $\Gamma$, and consider a function $\psi$ such that
$\psi(0)<0$, $\psi(T)>0$ and $\hat{I}_{T,x}(\psi)<+\infty$ (such a
function can be easily obtained by adapting the construction of the
function $\psi$ in the proof of Proposition~\ref{prop:I-not-lsc}
below).  Obviously, $\tilde{I}_{T,x}(\psi)=+\infty$, giving the
required counter-example.

Therefore, our upper bound is more precise than the one obtained by
classical general methods. This also explains why we have to use a
method based on a precise study of the paths of $X^{\varepsilon,x}$ to
obtain our result.

\paragraph{Proof of Proposition~\ref{prop:I-not-lsc}}
Take $y$ as in Proposition~\ref{prop:I-not-lsc}. By translation, we
can suppose that $y=0$. Then, Proposition~\ref{prop:d>1} implies that
there exists a neighborhood ${\cal N}_0$ of $0$ and a constant $a_0>0$
such that for all $s\in\mathbb{R}^d$ and $x\in{\cal N}_0$,
$s^*a(x)s\geq a_0\|x\|\|s\|^2$, \emph{i.e.} each eigenvalue of $a(x)$
is greater than $a_0\|x\|$. Therefore, for all $s\in\mathbb{R}^d$ and
$x\in{\cal N}_0$,
\begin{equation}
  \label{eq:I-not-lsc}
  s^*a^{-1}(x)s\leq\frac{\|s\|^2}{a_0\|x\|}.
\end{equation}

Take $x_0\in\RR^d\setminus\Gamma$ such that the segment $(0,x_0]$ is
included in $(\RR^d\setminus\Gamma)\cup{\cal N}_0$, and define for
$0\leq t\leq T$
\begin{equation*}
  \psi(t)=\left(1-\frac{2t}{T}\right)^2x_0,
\end{equation*}
and for all $n\geq 1$
\begin{equation*}
  \psi_n(t)=\left\{
    \begin{array}{ll}
      \psi(t) & \mbox{if}\quad t\in\displaystyle{\left[0,\frac{T}{2}-
          \frac{1}{n}\right]\cup\left[\frac{T}{2}+\frac{1}{n},T\right]} \\
      \displaystyle{\psi\left(\frac{T}{2}-\frac{1}{n}\right)}
      & \mbox{otherwise.}
    \end{array}\right.
\end{equation*}
Since $\psi(T/2-1/n)=\psi(T/2+1/n)$, $\psi_n$ is continuous and
piecewise differentiable.  Note that $\psi(t)$ and $\psi_n(t)$ belong
to $[0,x_0]$ for all $t\in[0,T]$, that $\psi(t)\not\in\Gamma$ except
if $t=T/2$, and that $\psi_n(t)\not\in\Gamma$ for any $t\in[0,T]$.
Therefore, $I_{T,x_0}(\psi)=\infty$, and $I_{T,x_0}(\psi_n)<\infty$.

It follows from~(\ref{eq:I-not-lsc}) and from the fact that $b$ is
$K$-Lipschitz that
\begin{equation}
  \label{eq:calc-prop-4.1}
  \begin{aligned}
    \hat{I}_{T,x_0}(\psi) & \leq\frac{1}{2a_0}\int_0^T
    \frac{\|(1-2t/T)2x_0/T+b(\psi(t))\|^2}{\|\psi(t)\|}dt \\
    & \leq\frac{1}{2a_0}\int_0^T\frac{2(1-2t/T)^24\|x_0\|^2/T^2
      +2K^2\|\psi(t)\|^2}{\|\psi(t)\|}dt \\
    & \leq\frac{1}{2a_0}\int_0^T\left(\frac{8}{T^2}\|x_0\|
      +2K^2\|\psi(t)\|\right)dt<\infty.
  \end{aligned}
\end{equation}
Now, for all $n\geq 1$,
\begin{align*}
  I_{T,x_0}(\psi_n) & \leq\hat{I}_{T,x_0}(\psi)+\frac{1}{2a_0}
  \int_{T/2-1/n}^{T/2+1/n}\frac{\|b(\psi_n(t))\|^2}{\|\psi_n(t)\|}dt
  \\ & \leq\hat{I}_{T,x_0}(\psi)+\frac{1}{2a_0}
  \int_{T/2-1/n}^{T/2+1/n}K^2\|\psi_n(t)\|dt,
\end{align*}
which is uniformly bounded in $n$. Hence $\limsup
I_{T,x_0}(\psi_n)<+\infty=I_{T,x_0}(\psi)$.

Let us extend this result to an arbitrary $x\not\in\Gamma$. Since the
points of $\Gamma$ are isolated in $\RR^d$, there exists $\alpha>0$
and $\phi\in{\cal C}^1([0,T],\Gamma_{\alpha})$ such that $\phi(0)=x$
and $\phi(T)=x_0$. Since $a$ is uniformly non-degenerate on
$\Gamma_{\alpha}$, $I_{T,x}(\phi)<\infty$. Therefore, it suffices to
concatenate $\phi$ and $\psi$ to obtain a function $\tilde{\psi}$
defined on $[0,2T]$ such that
$\limsup\tilde{I}_{2T,x}(\tilde{\psi})<I_{2T,x}(\tilde{\psi})$. Since
this can be done for all $T>0$, this ends the proof of
Proposition~\ref{prop:I-not-lsc}.\hfill$\Box$

\subsection{Proof of Theorem~\ref{thm:true-LDP}}
\label{sec:pf-true-LDP}

We first give some notation used throughout the proof.
\begin{itemize}
\item ${\cal C}_x(I,E) (resp.\ {\cal C}^{ac}_x(I,E),\ {\cal
    C}^1_x(I,E)$) is the set of continuous functions from
  $I\subset\mathbb{R}_+$ to $E\subset\mathbb{R}^d$ (resp.\ absolutely
  continuous, resp.\ ${\cal C}^1$) with value $x$ at $0$, endowed with
  the $L^\infty$ norm.
\item For $\varphi\in{\cal C}([0,T],\mathbb{R}^d)$ and $0\leq a<b\leq
  T$, define
  \begin{equation}
    \label{eq:def-|.|a,b}
    \|\varphi\|_{a,b}=\sup_{a\leq t\leq b}\|\varphi(t)\|,
  \end{equation}
  and
  \begin{equation}
    \label{eq:def-B}
    B_b(\varphi,\delta)=\{\tilde{\varphi}\in{\cal
      C}([0,T],\mathbb{R}^d):\|\tilde{\varphi}-\varphi\|_{0,b}\leq\delta\}.
  \end{equation}
  When $a=0$ and $b=T$, $\|\cdot\|_{0,T}$ is the usual $L^\infty$ norm
  in ${\cal C}([0,T],\mathbb{R}^d)$, and $B_T(\varphi,\delta)$ is the
  usual closed ball centered at $\varphi$ with radius $\delta$ for
  this norm, also simply denoted $B(\varphi,\delta)$. 
\end{itemize}
\bigskip

We are actually going to prove the following result.

\begin{thm}
  \label{thm:LDP}
  Assume the conditions of Theorem~\ref{thm:true-LDP}. Then, for any
  $x\in\RR^d$ and any open subset $O$ of ${\cal C}([0,T],\RR^d)$,
  \begin{equation}
    \label{eq:LB}
    \liminf_{\varepsilon\rightarrow
      0,y\rightarrow x}\varepsilon\ln\mathbb{P}(X^{\varepsilon,y}\in
    O)\geq -\inf_{\psi\in O}I_{T,x}(\psi),
  \end{equation}
  and for any $x\not\in\Gamma$ and any closed subset $C$ of ${\cal
    C}([0,T],\RR^d)$ such that ${\cal
    C}^1_x([0,T],\RR^d\setminus\Gamma)$ is dense in $C\cap{\cal
    C}_x([0,T],\RR^d)$,
  \begin{equation}
    \label{eq:UB}
    \limsup_{\varepsilon\rightarrow
      0,y\rightarrow x} \varepsilon\ln\mathbb{P}(X^{\varepsilon,y}\in
    C)\leq -\inf_{\psi\in C}I_{T,x}(\psi).
  \end{equation}
\end{thm}

This is an incomplete LDP involving the non-lower semicontinuous rate
function $I_{T,x}$. From this can be deduced the LDP involving the
rate function $\tilde{I}_{T,x}$ (Theorem~\ref{thm:true-LDP}) as
follows.

First, by definition of $\tilde{I}_{T,x}$, for any open $O\subset{\cal
  C}([0,T],\RR^d)$,
$$
\inf_{\psi\in O}I_{T,x}(\psi)=\inf_{\psi\in O}\tilde{I}_{T,x}(\psi).
$$
Therefore,~(\ref{eq:true-LB}) is immediate.

Moreover, $\tilde{I}_{T,x}\leq I_{T,x}$, so~(\ref{eq:true-UB})
obviously holds for the same closed sets as in
Theorem~\ref{thm:LDP}. Now, let $K$ be any compact subset of ${\cal
  C}([0,T],\RR^d)$. Since $\tilde{I}_{T,x}$ is lower semicontinuous,
for any $\eta>0$, there exists $\alpha>0$ such that
$$
\inf_{\psi\in K}\tilde{I}_{T,x}(\psi)\leq \inf_{\psi\in K_\alpha}
\tilde{I}_{T,x}(\psi) +\eta,
$$
where
$$
K_\alpha=\bigcup_{\psi\in K}B(\psi,\alpha).
$$
Indeed, if this would fail, there would exist $\eta>0$ and two
sequences $(\psi_n)_{n\geq 1}$ and $(\tilde{\psi}_n)_{n\geq 1}$ such
that $\tilde{\psi}_n\in K$, $\|\psi_n-\tilde{\psi}_n\|_{0,T}\leq 1/n$
and $\tilde{I}_{T,x}(\psi_n)\leq \tilde{I}_{T,x}(\tilde{\psi}_n)-\eta$
for all $n\geq 1$. Since $K$ is compact, we could then extract a
subsequence $(\tilde{\psi}_{i_n})$ of $(\tilde{\psi}_n)$ converging to
some $\tilde{\psi}\in K$. Since $\tilde{I}_{T,x}$ is lower
semicontinuous, this would imply that
$$
\tilde{I}_{T,x}(\tilde{\psi})
\leq\liminf_{n\rightarrow+\infty}\tilde{I}_{T,x}(\psi_{i_n}) \leq\inf_{\psi\in
  K}\tilde{I}_{T,x}(\psi)-\eta,
$$
which is a contradiction.

Now, let $\psi_1,\ldots,\psi_n$ be such that
$$
\tilde{K}_\alpha=\bigcup_{i=1}^nB(\psi_i,\alpha)\supset K.
$$
Since $\tilde{K}_\alpha\subset K_\alpha$,
$$
\inf_{\psi\in K}\tilde{I}_{T,x}(\psi)\leq \inf_{\psi\in \tilde{K}_\alpha}
\tilde{I}_{T,x}(\psi) +\eta.
$$
Moreover, the points of $\Gamma$ are isolated, and thus any point of
the interior of $\tilde{K}_\alpha$ is obviously limit of elements of
$\tilde{K}_\alpha\cap{\cal C}^1([0,T],\RR^d\setminus\Gamma)$. Since
$\tilde{K}_\alpha$ is the closure of its interior, any point of
$\partial\tilde{K}_\alpha$ is also limit of elements of
$\tilde{K}_\alpha\cap{\cal C}^1([0,T],\RR^d\setminus\Gamma)$ by a
diagonal procedure. Moreover, $\tilde{K}_\alpha$ is closed. Therefore,
one can apply~(\ref{eq:UB}) to this set:
\begin{align*}
  \limsup_{\varepsilon\rightarrow 0,\:y\rightarrow
    x}\varepsilon\ln\PP(X^{\varepsilon,y}\in K) & \leq
  \limsup_{\varepsilon\rightarrow 0,\:y\rightarrow
    x}\varepsilon\ln\PP(X^{\varepsilon,y}\in\tilde{K}_\alpha)\leq
  -\inf_{\psi\in\tilde{K}_\alpha}I_{T,x}(\psi) \\ & \leq
  -\inf_{\psi\in\tilde{K}_\alpha}\tilde{I}_{T,x}(\psi)\leq
  -\inf_{\psi\in K}\tilde{I}_{T,x}(\psi)+\eta.
\end{align*}
Since this holds for all $\eta>0$,~(\ref{eq:true-UB}) is proved for
compact sets.

The extension to any closed sets is classically deduced from the
following uniform exponential tightness estimate.
\begin{lemma}
  \label{lem:expo-tight}
  For any $k>0$ and $y\in\RR^d$, define the compact set
  \begin{equation}
    \label{eq:def_Kk}
    K^y_k=\left\{\psi\in{\cal C}_y([0,T],\RR^d):\forall l\geq k,\
      \omega\left(\psi,\frac{1}{l^3}\right)
      \leq\frac{1}{l}\right\},
  \end{equation}
  where $\omega(\psi,\delta)=\sup_{|t-s|\leq\delta}\|\psi(t)-
  \psi(s)\|$. Then, there exists $k_0$ and $\varepsilon_0$, such that
  for all $y\in\RR^d$, $k\geq k_0$ and
  $\varepsilon\leq\varepsilon_0$,
  \begin{equation}
    \label{eq:expo-tight}
    \varepsilon\ln\mathbb{P}(X^{\varepsilon,y}\not\in
    K^y_k)\leq -\frac{k}{64d\Sigma^2},
  \end{equation}
  where $\Sigma:=\sup_{x\in\RR^d}\|\sigma(x)\|$.
\end{lemma}

Then, taking any closed $C\in{\cal C}([0,T],\RR^d)$ and choosing $k$ large
enough,
\begin{align}
  \limsup_{\varepsilon\rightarrow 0,y\rightarrow x}\varepsilon\ln
  \mathbb{P}(X^{\varepsilon,y}\in C) &
  \leq\sup\left\{\limsup_{\varepsilon\rightarrow 0,y\rightarrow
      x}\varepsilon\ln\mathbb{P}(X^{\varepsilon,y}\in
    C\cap K^y_k),\right. \notag \\
  & \phantom{\leq\sup}\quad\left.\limsup_{\varepsilon\rightarrow 0,y\rightarrow
      x}\varepsilon\ln\mathbb{P}(X^{\varepsilon,y}\not\in K^y_k)
  \right\} \notag \\ & \leq -\inf_{\psi\in C}I_{T,x}(\psi),
  \label{eq:compact->closed}
\end{align}
ending the proof of Theorem~\ref{thm:true-LDP}.\hfill$\Box$
\bigskip

The proof of Lemma~\ref{lem:expo-tight} makes use of the following
classical exponential inequality for stochastic integrals, of which
the proof is omitted. This result will be also used in the proof of
Theorem~\ref{thm:LDP} below. Let ${\cal M}_{d,d}$ denote the set of
real $d\times d$ matrices.
\begin{lemma}
  \label{lem:Stroock}
  Let $Y_t$ be a ${\cal F}_t$-martingale with values in $\mathbb{R}^d$
  on a filtered probability space $(\Omega,{\cal F},{\cal
    F}_t,\mathbb{P})$, and suppose that its quadratic covariation
  process $\langle Y\rangle_t$ satisfies $\sup_{t\leq T}\|\langle
  Y\rangle_t\|\leq A$. Let $\tau$ be a ${\cal F}_t$ stopping time, and
  let $Z:\mathbb{R}_+\times\Omega\rightarrow{\cal M}_{d,d}$ be a
  progressively measurable process such that
  $\sup_{t\leq\tau}\|Z_t\|\leq B$. Then for any $R>0$,
  \begin{equation*}
    \mathbb{P}\left(\sup_{t\leq
    T}\left\|\int_0^{t\wedge\tau}Z_sdY_s\right\|\geq
    R\right)\leq 2d\exp\left(-\frac{R^2}{2dTAB^2}\right).
  \end{equation*}
\end{lemma}

\paragraph{Proof of Lemma~\ref{lem:expo-tight}}
It follows from~(\ref{eq:SDE-LDP}) that, for any $y\in\RR^d$, $s>0$
and $t\in[0,T]$,
\begin{equation*}
  \|X^{\varepsilon,y}_{t+s}-X^{\varepsilon,y}_t\|\leq
  Cs+\sqrt{\varepsilon}\left\|\int_t^{t+s}
    \sigma(X^{\varepsilon,y}_u)dW_u\right\|.
\end{equation*}
Fix $h>0$ and $R\geq Ch$. Applying Lemma~\ref{lem:Stroock}, we have
\begin{equation*}
  \mathbb{P}\left(\sup_{0\leq s\leq
  h}\|X^{\varepsilon,y}_{t+s}-X^{\varepsilon,y}_t\|\geq
  R\right)\leq 2d\exp\left(-\frac{(R-Ch)^2}{2dh\varepsilon\Sigma^2}\right).
\end{equation*}
Writing this for $t=ih$ for $0\leq i<T/h$, we deduce that
\begin{equation}
  \label{eq:expo-tight-1}
  \mathbb{P}\left(\omega(X^{\varepsilon},h)\geq
    2R\right)\leq 2d\left(\frac{T}{h}+1\right)\exp\left(
  -\frac{(R-Ch)^2}{2d\varepsilon\Sigma^2h}\right).
\end{equation}
For any $l\geq 1$, set $R_l=1/2l$ and $h_l=1/l^3$. Then, for
sufficiently large $l$, $R_l\geq Ch_l$ and
\begin{equation}
  \label{eq:expo-tight-2}
  \frac{(R_l-Ch_l)^2}{2d\varepsilon\Sigma^2h_l}
  =\frac{(\sqrt{l}-2C/l^{3/2})^2}{8d\varepsilon\Sigma^2}
  \geq\frac{(\sqrt{l}/2)^2}{8d\varepsilon\Sigma^2}
  =\frac{l}{32d\varepsilon\Sigma^2}.
\end{equation}
Observing that
\begin{equation*}
    K^y_k=\{\psi\in{\cal C}_y([0,T],{\cal X}):\forall l\geq k,
    \omega\left(\psi,h_l\right)\leq 2R_l\},
\end{equation*}
inequality~(\ref{eq:expo-tight}) easily follows
from~(\ref{eq:expo-tight-1}) and~(\ref{eq:expo-tight-2}).\hfill$\Box$

\subsection{Proof of Theorem~\ref{thm:LDP}}
\label{sec:pf-LDP}

The proof of Theorem~\ref{thm:LDP} makes use of the function $I_{T,x}$
and of the (good) rate function of Schilder's theorem (LDP for
Brownian motion)
\begin{equation*}
  J_{T}(\varphi)=\left\{
    \begin{array}{ll}
      \displaystyle{\frac{1}{2}\int_0^T\|\dot{\varphi}(t)\|^2dt} &
      \mbox{if}\quad\varphi\in{\cal C}^{ac}_0([0,T],\mathbb{R}^{d}) \\
      +\infty & \mbox{otherwise.}
    \end{array} \right.
\end{equation*}

First, we need to construct the function $S$ ``mapping'' Brownian
paths to the paths of $X^{\varepsilon}$. For any $\varphi\in{\cal
  C}_0^{ac}([0,T],\mathbb{R}^{d})$, let $S(\varphi)$ be the solution
on $[0,T]$ to
\begin{equation}
  \label{eq:def-S}
  S(\varphi)(t)=x+\int_0^tb(S(\varphi)(s))ds+
  \int_0^t\sigma(S(\varphi)(s))\dot{\varphi}(s)ds,
\end{equation}
obtained as follows: by Proposition~\ref{prop:a-b-btilde}~(i)
and~(iii), $b$ and $\sigma$ are bounded and locally Lipschitz on
$\RR^d\setminus\Gamma$. Therefore, Cauchy-Lipschitz's theorem implies
local existence and uniqueness in $\RR^d\setminus\Gamma$ of a solution
to $\dot{y}=b(y)+\sigma(y)\dot{\varphi}$. This defines properly
$S(\varphi)$ until the time $t_{S(\varphi)}$ where it reaches
$\Gamma$. In the case where $t_{S(\varphi)}<T$, set
$S(\varphi)(t)=S(\varphi)({t_{S(\varphi)}})$ for $t_{S(\varphi)}\leq
t\leq T$. The function $S(\varphi)$ obtained this way is a solution
to~(\ref{eq:def-S}) on $[0,T]$ and belongs to $\tilde{{\cal
    C}}^{ac}_x([0,T],\RR^d)$.  
\medskip

The proof of Theorem~\ref{thm:LDP} is based on the following two
lemmas. The first one gives a precise sense to the fact that the
function $S$ maps the paths of $\sqrt{\varepsilon}W$ to the paths of
$X^{\varepsilon,x}$. The second one gives the relation between $S$,
$I_{T,x}$ and $J_T$. Their proof is postponed after the proof of the
theorem.

\begin{lembr}
  \label{lem:approx}
  \begin{description}
  \item[\textmd{(i)}] Fix $\varphi\in{\cal
      C}^{ac}_0([0,T],\mathbb{R}^d)$ such that $\psi:=S(\varphi)$
    takes no value in $\Gamma$ and such that $J_T(\varphi)<+\infty$.
    Then, for all $\eta>0$ and $R>0$, there exists $\delta>0$ such
    that
  \begin{equation}
    \label{eq:lem-approx-(i)}
    \limsup_{\varepsilon\rightarrow 0,y\rightarrow x}\varepsilon\ln
    \mathbb{P}\left(\|X^{\varepsilon,y}
      -S(\varphi)\|_{0,T}\geq\eta,
      \|\sqrt{\varepsilon}W-\varphi\|_{0,T}
      \leq\delta\right)\leq -R.
  \end{equation}
  \item[\textmd{(ii)}] Fix $\tilde{\varphi}\in{\cal
      C}^{ac}_0([0,T],\mathbb{R}^d)$ such that
    $\psi(t):=S(\tilde{\varphi})(t)\in\Gamma$ for some $t\in[0,T]$.
    Define $\varphi(t)=\tilde{\varphi}(t)$ for $t<t_{\psi}$ and
    $\varphi(t)=\tilde{\varphi}(t_{\psi})$ for $t_{\psi}\leq t\leq T$.
    Then $S(\varphi)=S(\tilde{\varphi})=\psi$.  Suppose that
    $J_T(\varphi)<+\infty$. Then, for all $\eta>0$ and $R>0$, there
    exists $\delta>0$ such that~(\ref{eq:lem-approx-(i)}) holds.
  \item[\textmd{(iii)}] With the same $\varphi$ as in~(i),
    for all $\delta>0$ and $R>0$, there exists $\eta>0$ such that
  \begin{equation}
    \label{eq:lem-approx-(ii)}
    \limsup_{\varepsilon\rightarrow 0,y\rightarrow x}\varepsilon\ln
    \mathbb{P}\left(\|X^{\varepsilon,y}
      -S(\varphi)\|_{0,T}\leq\eta,
      \|\sqrt{\varepsilon}W-\varphi\|_{0,T}
      \geq\delta\right)\leq -R.
  \end{equation}
  \end{description}
\end{lembr}

\begin{lembr}
  \label{lem:S-I-Itilde}
  \begin{description}
  \item[\textmd{(i)}] For all $\psi\in{\cal C}_x([0,T],\RR^d)$,
    $$
    I_{T,x}(\psi)=\inf\{ J_{T}(\varphi),S(\varphi)=\psi\}
    $$
    and when $I_{T,x}(\psi)<+\infty$, there is a unique $\varphi\in{\cal
      C}^{ac}_0([0,T],\mathbb{R}^d)$ that realizes this infimum, and
    this function is constant after $t_{\psi}$.
  \item[\textmd{(ii)}] ${\cal C}^1([0,T],\RR^d\setminus\Gamma)$ is
    dense in $S(\{J_T<\infty\})$.
  \end{description}
\end{lembr}

In~\cite{azencott-80,doss-priouret-83}, $b^{\varepsilon}$ and $\sigma$
are assumed Lipschitz, and thus Point~(i) of Lemma~\ref{lem:approx}
can be proved for all $\varphi\in{\cal C}^{ac}_0([0,T],\RR^d)$, which
is enough to conclude. In our case, because of the bad regularity of
the coefficients of the SDE, we cannot prove~(i) for all
$\varphi\in{\cal C}^{ac}_0([0,T],\RR^d)$. As a consequence, we are
only able to obtain the large deviations lower bound from
Lemma~\ref{lem:approx}~(i) and~(ii). In order to prove the large
deviations upper bound, we use an original method based on
Lemma~\ref{lem:approx}~(iii).

Lemma~\ref{lem:S-I-Itilde} is an extension to our case of very similar
lemmas in~\cite{azencott-80,doss-priouret-83}.

\paragraph{Proof of Theorem~\ref{thm:LDP}: lower bound}
It is well-known that the lower bound~(\ref{eq:LB}) for any open set
$O$ is equivalent to the fact that, for all $\psi\in{\cal
  C}_x([0,T],\RR^d)$ and $\eta>0$,
\begin{equation}
  \label{eq:LB-1}
  \liminf_{\varepsilon\rightarrow 0,y\rightarrow x}\varepsilon\ln
  \mathbb{P}(\|X^{\varepsilon,y}-\psi\|_{0,T}
  \leq\eta)\geq -I_{T,x}(\psi).
\end{equation}
Fix $\psi$ and $\eta$ as above, and assume that
$I_{T,x}(\psi)<+\infty$ (otherwise, there is nothing to prove). By
Lemma~\ref{lem:S-I-Itilde}~(i), there is a unique $\varphi\in{\cal
  C}^{ac}_0([0,T],\mathbb{R}^d)$ such that $S(\varphi)=\psi$ and
$u:=J_{T}(\varphi)=I_{T,x}(\psi)$. Choose $R>u$.  If the image of
$\psi$ has empty intersection with $\Gamma$, apply
Lemma~\ref{lem:approx}~(i). Otherwise, apply
Lemma~\ref{lem:approx}~(ii). In both cases, there exists $\delta>0$
such that
\begin{equation*}
  \limsup_{\varepsilon\rightarrow 0,y\rightarrow x}\varepsilon\ln
    \mathbb{P}\left(\|X^{\varepsilon,y}-\psi\|_{0,T}
      \geq\eta, \|\sqrt{\varepsilon}W-\varphi\|_{0,T}
      \leq\delta\right)\leq -R.
\end{equation*}
Since
\begin{align*}
  \mathbb{P}(\|\sqrt{\varepsilon}
    W-\varphi\|_{0,T}\leq\delta)
  & \leq\mathbb{P}(\|X^{\varepsilon,y}
  -\psi\|_{0,T}<\eta) \\
  & +\mathbb{P}(\|X^{\varepsilon,y}-\psi\|_{0,T}
  \geq\eta, \|\sqrt{\varepsilon}W-\varphi\|_{0,T}
  \leq\delta),
\end{align*}
we deduce from Schilder's theorem that
\begin{align*}
  -u=-J_{T}(\varphi) & \leq-\inf\{J_{T}(\tilde{\varphi}),
  \tilde{\varphi}\in B_T(\varphi,\delta)\} \\
  & \leq\liminf_{\varepsilon\rightarrow 0,y\rightarrow x}\varepsilon
  \ln\mathbb{P}(\|\sqrt{\varepsilon}
  W-\varphi\|_{0,T}<\delta) \\
  & \leq\sup\left\{\liminf_{\varepsilon\rightarrow 0,y\rightarrow x}
    \varepsilon\ln\mathbb{P}(\|X^{\varepsilon,y}
    -\psi\|_{0,T}<\eta),\right. \\
  & \phantom{\sup}\left.\quad\liminf_{\varepsilon\rightarrow 0,y\rightarrow x}
    \varepsilon\ln\mathbb{P}(\|X^{\varepsilon,y}
    -\psi\|_{0,T}\geq\eta,\|\sqrt{\varepsilon}W
    -\varphi\|_{0,T}\leq\delta)\right\} \\
  & \leq\sup\left\{\liminf_{\varepsilon\rightarrow 0,y\rightarrow x}
    \varepsilon\ln\mathbb{P}(\|X^{\varepsilon,y}
    -\psi\|_{0,T}<\eta),-R\right\},
\end{align*}
and since $R>u$,~(\ref{eq:LB-1}) is established.\hfill$\Box$

\paragraph{Proof of Theorem~\ref{thm:LDP}: upper bound}
We first prove~(\ref{eq:UB}) for particular compact sets: let $K$ be a
non-empty compact set of ${\cal C}([0,T],\RR^d)$ such that
$S(\{J_T<+\infty\})$ is dense in $K_x$, where $K_x:=K\cap{\cal
  C}_x([0,T],\RR^d)$. By Lemma~\ref{lem:S-I-Itilde}~(i),
$S(\{J_T<+\infty\})=\{I_{T,x}<+\infty\}$, and so $u:=\inf\{
I_{T,x}(\psi),\psi\in K\}<+\infty$.

Fix $\rho>0$. For any $\psi\in K\cap S(\{J_T<+\infty\})$, by
Lemma~\ref{lem:S-I-Itilde}~(i), there exists a unique $\varphi\in{\cal
  C}_0^{ac}([0,T], \mathbb{R}^d)$ constant after $t_\psi$ such that
$S(\varphi)=\psi$ and $I_{T,x}(\psi)=J_{T}(\varphi)<\infty$. We intend
to use Lemma~\ref{lem:approx}~(iii), which holds only if $\psi$ takes
no value in $\Gamma$. So we have to introduce $\alpha_{\psi}>0$ such
that
\begin{equation*}
  \frac{1}{2}\int_{t_{\psi}-\alpha_{\psi}}^{t_{\psi}}
  \|\dot{\varphi}_s\|^2ds<\frac{\rho}{2},
\end{equation*}
so that $J_T(\varphi)\leq
J_{t_{\psi}-\alpha_{\psi}}(\varphi)+\rho/2$. Since
$J_{t_{\psi}-\alpha_{\psi}}$ is lower semicontinuous, there exists
$\delta_{\psi}>0$ such that
\begin{equation}
  \label{eq:ineq-UB-1}
  \forall\tilde{\varphi}\in B_{t_{\psi}-\alpha_{\psi}}
  (\varphi,\delta_{\psi}),\quad
  J_{t_{\psi}-\alpha_{\psi}}(\tilde{\varphi})\geq
  J_{t_{\psi}-\alpha_{\psi}}(\varphi)-\frac{\rho}{2}\geq
  J_T(\varphi)-\rho,
\end{equation}
where $B_t(\varphi,\delta)$ has been defined in~(\ref{eq:def-B}).

Applying Lemma~\ref{lem:approx}~(ii) to $\psi$ with
$T=t_{\psi}-\alpha_{\psi}$, $\delta=\delta_{\psi}$ and $R>u$, there
exists $\eta_{\psi}>0$ such that
\begin{equation}
 \label{eq:inequ-UB-2}
  \limsup_{\varepsilon\rightarrow 0,y\rightarrow x}\varepsilon\ln
    \mathbb{P}\left(\|X^{\varepsilon,y}
      -\psi\|_{0,t_{\psi}-\alpha_{\psi}}\leq\eta_{\psi},
      \|\sqrt{\varepsilon}W-\varphi\|_{0,t_{\psi}-
    \alpha_{\psi}}\geq\delta_{\psi}\right)\leq -R.
\end{equation}

Since we have assumed that $K_x\cap S(\{ J_T<+\infty\})$ is dense in
$K_x$, 
\begin{equation*}
  K_x\subset\bigcup_{\psi\in
  K_x\cap S(\{ J_T<+\infty\})}B_T(\psi,\eta_{\psi}).
\end{equation*}
Since $K_x$ is compact, there exists a finite number of functions
$\psi_1,\ldots,\psi_n$ in $K_x\cap S(\{ J_T<+\infty\} )$ such that
\begin{equation*}
  K_x\subset\bigcup_{i=1}^nB_T(\psi_i,\eta_i),
\end{equation*}
where we wrote $\eta_i$ instead of $\eta_{\psi_i}$. Since $K$ is
compact, there exists a neighborhood ${\cal N}_x$ of $x$ such that
\begin{equation*}
  K_{{\cal N}_x}\subset\bigcup_{i=1}^nB_T(\psi_i,\eta_i),
\end{equation*}
where $K_{{\cal N}_x}=\{\psi\in K:\psi(0)\in{\cal N}_x\}$. 

Now, define
\begin{equation*}
  U=\bigcup_{i=1}^nB_{t_i-\alpha_i}(\varphi_i,\delta_i),
\end{equation*}
where $t_i=t_{\psi_i}$, $\alpha_i=\alpha_{\psi_i}$ and
$\delta_i=\delta_{\psi_i}$, and where $\varphi_i$ is the function
satisfying $S(\varphi_i)=\psi_i$ and
$I_{T,x}(\psi_i)=J_{T}(\varphi_i)$. Then, for any $y\in{\cal N}_x$,
\begin{align*}
  \mathbb{P}(X^{\varepsilon,y}\in K) &
  \leq\mathbb{P}(\sqrt{\varepsilon}W\in
  U)+\mathbb{P}(\sqrt{\varepsilon}W\not\in
  U,\ X^{\varepsilon,y}\in K_{{\cal N}_x}) \\
  & \leq\sum_{i=1}^n\mathbb{P}(\sqrt{\varepsilon}
  W\in B_{t_i-\alpha_i}(\varphi_i,\delta_i)) \\
  & \quad+\sum_{i=1}^n\mathbb{P}(\|X^{\varepsilon,y}
  -\psi_i\|_{0,T}\leq\eta_i,\sqrt{\varepsilon}W\not\in U) \\
  & \leq\sum_{i=1}^n\mathbb{P}(\|\sqrt{\varepsilon}
  W-\varphi_i\|_{0,t_{i}-\alpha_i}<\delta_i) \\
  & \quad+\sum_{i=1}^n\mathbb{P}(\|X^{\varepsilon,y}
  -\psi_i\|_{0,t_{i}-\alpha_i}\leq\eta_i,\|\sqrt{\varepsilon}
  W-\varphi_i\|_{0,t_{i}-\alpha_i}\geq\delta_i).
\end{align*}
Since by Schilder's Theorem and~(\ref{eq:ineq-UB-1})
\begin{equation*}
  \limsup_{\varepsilon\rightarrow 0,y\rightarrow x}
    \varepsilon\ln\mathbb{P}(\|\sqrt{\varepsilon}
    W-\varphi_i\|_{0,t_i-\alpha_i}\leq\delta_i)
  \leq-\inf_{\varphi\in B_{t_i-\alpha_i}(\varphi_i,\delta_i)}
  J_{t_i-\alpha_i}(\varphi)\leq-J_T(\varphi_i)+\rho,
\end{equation*}
we finally deduce from~(\ref{eq:inequ-UB-2}) that
\begin{align*}
  \limsup_{\varepsilon\rightarrow 0,y\rightarrow x}\varepsilon\ln
  \mathbb{P}(X^{\varepsilon,y}\in K) & \leq\sup\left\{\sup_{1\leq i\leq
      n}(-J_{T}(\varphi_i)+\rho),
    -R\right\} \\
  & \leq\sup\left\{-\inf\{ I_{T,x}(\psi),\psi\in
    K\}+\rho,-R\right\}\leq-u+\rho.
\end{align*}
Since this holds for any $\rho>0$, the proof of~(\ref{eq:UB}) for the
set $K$ is completed.
\bigskip

Now, let $C$ be a closed subset of ${\cal C}([0,T],\RR^d)$ such that
${\cal C}^1_x([0,T],\RR^d\setminus\Gamma)$ is dense in $C\cap{\cal
  C}_x([0,T],\RR^d)$. Define the compact set
\begin{align*}
  K_k & =\{\psi\in{\cal C}([0,T],\RR^d):\|\psi(0)-x\|\leq
  1,\forall l\geq k,\omega(\psi,1/k^3)\leq 1/k\} \\
  & =\bigcup_{\|y-x\|\leq 1}K^y_k,
\end{align*}
where $K^y_k$ is defined in~(\ref{eq:def_Kk}). In order to apply the
previous upper bound for compact sets, we are going to construct a
compact set $\tilde{K}_k\supset K_k$ such that $S\{(J_T<\infty\})$ is
dense in $C\cap\tilde{K}_k\cap{\cal C}_x([0,T],\RR^d)$. This will be
enough to conclude since, by
Lemma~\ref{lem:expo-tight},
\begin{equation}
  \label{eq:expo-tight-UB-1}
  \limsup_{\varepsilon\rightarrow 0,y\rightarrow x}\varepsilon\ln
  \mathbb{P}(X^{\varepsilon,y}\not\in \tilde{K}_k)\leq
  \limsup_{\varepsilon\rightarrow 0,y\rightarrow x}\varepsilon\ln
  \mathbb{P}(X^{\varepsilon,y}\not\in K_k)\leq -k/64d\Sigma^2,
\end{equation}
so that the upper bound~(\ref{eq:UB}) will be proved as
in~(\ref{eq:compact->closed}).

The set $\tilde{K}_k$ can be constructed as follows. The set $C\cap
K_k\cap{\cal C}_x([0,T],\RR^d)$ is compact, so it is separable. Let
$(\psi_n)_{n\geq o}$ be a sequence of functions dense in this set.
For all $n\geq 0$, $\psi_n\in C$, so, by assumption, there exists a
sequence $(\psi_{n,p})_{p\geq 0}$ in $C\cap{\cal
  C}^1_x([0,T],\RR^d\setminus\Gamma)$ converging to $\psi_n$, such
that $\|\psi_{n,p}-\psi_n\|_{0,T}\leq 2^{-p}$ for all $p\geq 0$. Let
us define
\begin{equation*}
  \tilde{K}_k=K_k\cup\left(\bigcup_{n\geq 0}\{\psi_{n,p}:p\geq n\}\right),
\end{equation*}
and let us prove that $\tilde{K}_k$ is compact. Let $(\phi_m)$ be a
sequence of $\tilde{K}_k$. Extracting a converging subsequence is
trivial, except in the case where $\{m:\phi_m\in K_k\}$ is finite, and
when for all $n\geq 0$, $\{m:\phi_m\in\{\psi_{n,p}:p\geq n\}\}$ is
finite. In this case, there exists two increasing sequences of
integers $(\alpha_m)$ and $(\beta_m)$ such that for all $m\geq 0$,
$\phi_{\alpha_m}\in\{\psi_{\beta_m,p}:p\geq\beta_m\}$. For all $m\geq
0$, $\psi_{\beta_m}$ belongs to the compact set $C\cap K_k\cap{\cal
  C}_x([0,T],\RR^d)$, so, extracting a subsequence from $(\beta_m)$,
we can assume that $\psi_{\beta_m}\rightarrow\psi\in C\cap
K_k\cap{\cal C}_x([0,T],\RR^d)$. Then
\begin{equation*}
  \|\phi_{\alpha_m}-\psi\|_{0,T}\leq
  2^{-\beta_m}+\|\psi_{\beta_m}-\psi\|\rightarrow 0
\end{equation*}
when $m\rightarrow\infty$. Hence $\tilde{K}_k$ is compact. Moreover,
$\tilde{K}_k$ has been constructed in such a way that ${\cal
  C}_x^1([0,T],\RR^d\setminus\Gamma)$ is dense in
$C\cap\tilde{K}_k\cap{\cal C}_x([0,T],\RR^d)$, as required. This ends
the proof of Theorem~\ref{thm:LDP}.\hfill$\Box$

\subsection{Proof of Lemmas~\ref{lem:approx} and~\ref{lem:S-I-Itilde}}
\label{sec:pf-lemmas}

\paragraph{Proof of Lemma~\ref{lem:approx}}
Let $\varphi$ be as in any point of Lemma~\ref{lem:approx}. We will
first restrict ourselves to the case $\varphi=0$ by means of
Girsanov's Theorem. Define on $(\Omega,{\cal F}_T)$ the
probability measure ${\mathbb{P}}^{\varepsilon,y}$ by
\begin{equation*}
  \frac{d{\mathbb{P}}^{\varepsilon,y}}{d\mathbb{P}}=
  \exp\left(\frac{1}{\sqrt{\varepsilon}}\int_0^T\dot{\varphi}_s
    dW_s-\frac{1}{2\varepsilon}\int_0^T\left\|\dot{\varphi}_s
    \right\|^2ds\right).
\end{equation*}
Since in all cases $J_T(\varphi)=1/2\int_0^T\|\dot{\varphi}_t\|^2dt
<+\infty$, by Novikov's criterion, Girsanov's Theorem is applicable
and implies that 
$$
\tilde{W}^{\varepsilon}_t:=W_t-
\frac{\varphi_t}{\sqrt{\varepsilon}}
$$
is a ${\mathbb{P}}^{\varepsilon,y}$-Brownian motion for $t\leq T$ and
that, ${\mathbb{P}}^{\varepsilon,y}$-a.s., for any $t\leq T$,
\begin{equation}
  \label{eq:Xtilde}
  {X}^{\varepsilon,y}_t=y+\int_0^t
  (b^{\varepsilon}({X}^{\varepsilon,y}_s)+
    \sigma({X}^{\varepsilon,y}_s)\dot{\varphi}_s)ds+
  \sqrt{\varepsilon}\int_0^t
  \sigma({X}^{\varepsilon,y}_s)d\tilde{W}^{\varepsilon}_s.
\end{equation}
Let 
\begin{align*}
  F^{\varepsilon,y} & =\{\|X^{\varepsilon,y}-S(\varphi)\|_{0,T}\geq\eta,\
  \|\sqrt{\varepsilon}W-\varphi\|_{0,T}\leq\delta\} \\ & =
  \{\|{X}^{\varepsilon,y}-S(\varphi)\|_{0,T}\geq\eta,\
  \|\sqrt{\varepsilon}\tilde{W}^{\varepsilon}\|_{0,T}\leq\delta\}.
\end{align*}
It follows from Cauchy-Schwartz's inequality that
\begin{equation}
  \label{eq:CS-Girsanov}
  \mathbb{P}(F^{\varepsilon,y})=\int
  \mathbf{1}_{F^{\varepsilon,y}}\frac{d\mathbb{P}}
  {d\mathbb{P}^{\varepsilon,y}}d{\mathbb{P}}^{\varepsilon,y}
  \leq\left({\mathbb{P}}^{\varepsilon,y}(F^{\varepsilon,y})
  \right)^{\frac{1}{2}}\left(\int\left(\frac{d\mathbb{P}}
      {d{\mathbb{P}^{\varepsilon,y}}}\right)^2
    d{\mathbb{P}}^{\varepsilon,y}\right)^{\frac{1}{2}}.
\end{equation}
Now,
\begin{equation*}
  \begin{aligned}
    \left(\frac{d\mathbb{P}}
      {d{\mathbb{P}}^{\varepsilon,y}}\right)^2 &
    =\exp\left(-\frac{2}{\sqrt{\varepsilon}}\int_0^T\dot{\varphi}_s
      d\tilde{W}_s^{\varepsilon}-\frac{1}{\varepsilon}
      \int_0^T\|\dot{\varphi}_s\|^2ds\right)\\
    & =\exp\left(\int_0^T\left(-\frac{2\dot{\varphi}_s}
        {\sqrt{\varepsilon}}\right)
      d\tilde{W}^{\varepsilon}_s-\frac{1}{2}\int_0^T\left\|
        \frac{2\dot{\varphi}_s}{\sqrt{\varepsilon}}\right\|^2ds\right) \\
    & \qquad\times\exp\left(\frac{1}{\varepsilon}
      \int_0^T\|\dot{\varphi}_s\|^2ds\right).
  \end{aligned}
\end{equation*}
The first term in the product of the right-hand side is a
${\mathbb{P}}^{\varepsilon,y}$-martingale (by Novikov's
criterion), and the second term is equal to
$\exp(2J_T(\varphi)/\varepsilon)$.  Therefore,~(\ref{eq:CS-Girsanov})
implies
\begin{equation*}
  \varepsilon\ln\mathbb{P}(F^{\varepsilon,y})\leq
  \frac{\varepsilon}{2}\ln{\mathbb{P}}^{\varepsilon,y}
    (F^{\varepsilon,y})+J_T(\varphi).
\end{equation*}
Therefore, Lemma~\ref{lem:approx} follows from the next
result.\hfill$\Box$

\begin{lemma}
  \label{lem:phi=0}
  The three points of Lemma~\ref{lem:approx} hold
  when~(\ref{eq:lem-approx-(i)}) and~(\ref{eq:lem-approx-(ii)}) are
  replaced respectively by
  \begin{equation}
    \label{eq:lem-phi=0-(i)}
    \limsup_{\varepsilon\rightarrow 0,y\rightarrow x}\varepsilon\ln
    {\mathbb{P}}^{\varepsilon,y}\left(\|{X}^{\varepsilon,y}-
      S(\varphi)\|_{0,T}
      \geq\eta,\|\sqrt{\varepsilon}\tilde{W}^{\varepsilon}\|_{0,T}
      \leq\delta\right)\leq-R
  \end{equation}
  and
  \begin{equation}
    \label{eq:lem-phi=0-(ii)}
    \limsup_{\varepsilon\rightarrow 0,y\rightarrow x}\varepsilon\ln
    {\mathbb{P}}^{\varepsilon,y}\left(\|{X}^{\varepsilon,y}-
      S(\varphi)\|_{0,T}
      \leq\eta,\|\sqrt{\varepsilon}\tilde{W}^{\varepsilon}\|_{0,T}
      \geq\delta\right)\leq -R.
  \end{equation}
\end{lemma}

Lemma~\ref{lem:phi=0} relies on the following lemma, of which the
proof is postponed after the proof of Lemma~\ref{lem:phi=0}.
\begin{lemma}
  \label{lem:prelim}
  With the previous notation, let $Y_t$ be a
  $\PP^{\varepsilon,y}$-martingale in $L^2$ such that $\sup_{t\leq
    T}\|\langle Y\rangle_t\|\leq A$, let $\tau$ be a stopping time,
  and let $\xi$ be a uniformly continuous bounded function on $\RR^d$.
  Then, for any $\eta>0$ and $R>0$, there exists $\delta>0$ and
  $\varepsilon_0>0$ both depending on $Y$ only through $A$ and both
  independent of $\tau$, such that for any $y\in\RR^d$ and
  $\varepsilon<\varepsilon_0$,
  \begin{equation}
    \label{eq:lem-prelim}
    \varepsilon\ln\mathbb{P}^{\varepsilon,y}\left(\left\|\sqrt{\varepsilon}
        \int_0^{\cdot\wedge\tau}\xi({X}^{\varepsilon,y}_s)
        dY_s\right\|_{0,T}\geq\eta,\ \|\sqrt{\varepsilon}Y\|_{0,T}
      \leq\delta\right)\leq -R.
  \end{equation}
\end{lemma}

\paragraph{Proof of Lemma~\ref{lem:phi=0}~(i)}
The function $\psi=S(\varphi)$ does not take any value in $\Gamma$ on
$[0,T]$, so there exists $\alpha>0$ such that $\forall t\in[0,T]$,
$\psi_t\in\Gamma_{\alpha}$. Suppose without loss of generality that
$\eta<\alpha/2$, and define for $y\in\RR^{d}$
\begin{equation*}
  \tau^{\varepsilon,y}=\inf\{ t:d({X}^{\varepsilon,y}_t,\Gamma)
  \leq\alpha/2\}\wedge T.
\end{equation*}
When $\tau^{\varepsilon,y}<T$, $\|{X}^{\varepsilon,y}_{
  \tau^{\varepsilon,y}}-S(\varphi)_{\tau^{\varepsilon,y}}\|\geq
d(S(\varphi)_{\tau^{\varepsilon,y}},\Gamma)-d({X}^{
  \varepsilon,y}_{\tau^{\varepsilon,y}},\Gamma)\geq\alpha/2>\eta$, so
\begin{equation*}
  \|{X}^{\varepsilon,y}-S(\varphi)\|_{0,T}
  \geq\eta\Rightarrow\|{X}^{\varepsilon,y}-
  S(\varphi)\|_{0,\tau^{\varepsilon,y}}\geq\eta.
\end{equation*}
Consequently,~(\ref{eq:lem-phi=0-(i)}) will be proved if we find
$\delta>0$ such that
\begin{equation*}
  \limsup_{\varepsilon\rightarrow 0,y\rightarrow x}\varepsilon
  \ln\mathbb{P}^{\varepsilon,y}(\|{X}^{\varepsilon,y}-S(\varphi)\|_{0,
    \tau^{\varepsilon,y}}\geq\eta,\ \|\sqrt{\varepsilon}W^\varepsilon\|_{0,T}
  \leq\delta)\leq-R.
\end{equation*}

Take $C$ such that $\sigma$ and $b$ are $C$-Lipschitz and $\tilde{b}$
is bounded by $C$ on $\Gamma_{\alpha/2}$.  It follows
from~(\ref{eq:Xtilde}) that, for $t\leq\tau^{\varepsilon,y}$,
\begin{align*}
  {\|{X}_t^{\varepsilon,y}-S(\varphi)_t\|} &
  \leq\sqrt{\varepsilon}\left\|\int_0^t\sigma
    ({X}_s^{\varepsilon,y})dW^\varepsilon_s\right\|
  +\varepsilon\int_0^t\|\tilde{b}({X}_s^{\varepsilon,y})\|ds+\|x-y\|
  \\ & \qquad+\int_0^t\|b({X}_s^{\varepsilon,y})-
    b(S(\varphi)_s)\|ds+\int_0^t\|\sigma({X}_s^{\varepsilon,y})-
    \sigma(S(\varphi)_s)\|\:\|\dot{\varphi}_s\|ds \\ & 
  \leq\sqrt{\varepsilon}\left\|\int_0^t
    \sigma({X}_s^{\varepsilon,y})dW^\varepsilon_s\right\|
  +\varepsilon CT+\|x-y\| +C\int_0^t(1+\|\dot{\varphi}_s\|)
  \|{X}_s^{\varepsilon,y}-S(\varphi)_s\|ds.
\end{align*}
Since $u:=\int_0^T\|\dot{\varphi}_s\|^2ds<+\infty$, by Gronwall's
lemma and the Cauchy-Schwartz's inequality, for
$t\leq\tau^{\varepsilon,y}$
\begin{equation*}
  {\|{X}_t^{\varepsilon,y}-S(\varphi)_t\|} 
  \leq\left(\sqrt{\varepsilon}\left\|\int_0^t
      \sigma({X}_s^{\varepsilon,y})dW^\varepsilon_s\right\|+\varepsilon
    CT+\|x-y\|\right)\exp\left(C\left(T+\sqrt{uT}\right)\right).
\end{equation*}
Therefore, it suffices to find $\delta>0$ such that
\begin{equation*}
  \limsup_{\varepsilon\rightarrow 0,y\rightarrow x}\varepsilon
  \ln\mathbb{P}^{\varepsilon,y}\left(\sqrt{\varepsilon}\left\|\int_0^t\sigma
      ({X}^{\varepsilon,y}_s)dW^\varepsilon_s\right\|_{0,\tau^{
        \varepsilon,y}}\geq\eta\beta,\ \sqrt{\varepsilon}\|W\|_{0,T}
    \leq\delta\right)\leq-R,
\end{equation*}
where $\beta=\exp[-C(T+\sqrt{uT})]/2$. This is an direct consequence
of Lemma~\ref{lem:prelim} with $Y=W^\varepsilon$, $A=1$, $\xi=\sigma$
and $\tau=\tau^{\varepsilon,y}$.\hfill$\Box$

\paragraph{Proof of Lemma~\ref{lem:phi=0}~(ii)}
In Lemma~\ref{lem:phi=0}~(ii), $\varphi$ is defined from
$\tilde{\varphi}$ by $\varphi_t=\tilde{\varphi}_t$ for $t\leq
t_{\psi}$, and $\varphi_t=\tilde{\varphi}_{t_{\psi}}$ otherwise, where
$\psi=S(\tilde{\varphi})=S(\varphi)$. By Cauchy-Schwartz's inequality,
$\int_0^{t_{\psi}}\|\dot{\varphi}_s\|ds\leq(2TJ_T(\varphi))^{1/2}<+\infty$,
so there exists $\rho>0$ small enough such that
\begin{equation}
  \label{eq:pf-lem2(iii)-1}
  \int_{t_{\psi}-\rho}^{t_{\psi}}\|\dot{\varphi}_s\|ds
  \leq\frac{\eta e^{-CT}}{8C},
\end{equation}
where $C$ is a constant bounding $b$, $\tilde{b}$ and $\sigma$, and
such that $b$ is $C$-Lipschitz.

Now, we have
\begin{equation*}
  \{\|{X}^{\varepsilon,y}-\psi\|_{0,T}\geq\eta,\
  \|\sqrt{\varepsilon}W^\varepsilon\|_{0,T}\leq\delta\}\subset
  D^{\varepsilon,y}\cup E^{\varepsilon,y},
\end{equation*}
where
\begin{gather*}
  D^{\varepsilon,y}=\left\{\|{X}^{\varepsilon,y}-
    \psi\|_{0,t_{\psi}-\rho}\leq\frac{\eta e^{-CT}}{4},\ 
    \|{X}^{\varepsilon,y}-\psi\|_{t_{\psi}-\rho,T}\geq\eta,
    \|\sqrt{\varepsilon}W^\varepsilon\|_{0,T}\leq\delta\right\} \\
  \mbox{and}\quad
  E^{\varepsilon,y}=\left\{\|{X}^{\varepsilon,y}-
    \psi\|_{0,t_{\psi}-\rho}\geq\frac{\eta e^{-CT}}{4},\ 
    \|\sqrt{\varepsilon}W^\varepsilon\|_{0,t_{\psi}-\rho}\leq\delta\right\}.
\end{gather*}
Part~(i) of Lemma~\ref{lem:phi=0} shows that
$\mathbb{P}^{\varepsilon,y}(E^{\varepsilon,y})$ has the required
exponential decay if $\delta$ is small enough. Let us estimate
$\PP^{\varepsilon,y}(D^{\varepsilon,y})$.

It follows from~(\ref{eq:Xtilde}) and from the fact that
$\dot{\varphi}_t=0$ for $t>t_{\psi}$ that, for any $t\geq
t_{\psi}-\rho$
\begin{align*}
  \|{X}^{\varepsilon,y}_t-\psi_t\|
  & \leq\|{X}^{\varepsilon,y}_{t_{\psi}-\rho}-
  \psi_{t_{\psi}-\rho}\|+\sqrt{\varepsilon}\left\|\int_{t_{\psi}-\rho}^t
    \sigma({X}^{\varepsilon,y}_s)dW^\varepsilon_s\right\| \\ &
  +C\int_{t_{\psi}-\rho}^t\|{X}^{\varepsilon,y}_s-\psi_s\|ds
  +\varepsilon CT+\int_{t_{\psi}-\rho}^{t_{\psi}\wedge
    t}\|\sigma({X}^{\varepsilon,y}_s)-
  \sigma(\psi_s)\|\:]\|\dot{\varphi}_s\|ds.
\end{align*}
On the event $D^{\varepsilon,y}$, the first term of the right-hand
side is smaller than $\eta e^{-CT}/4$, and, since $\sigma$ is bounded
by $C$, the last term is smaller than
$2C\int_{t_{\psi}-\rho}^{t_{\psi}}\|\dot{\varphi}\|ds\leq \eta
e^{-CT}/4$ by~(\ref{eq:pf-lem2(iii)-1}).  Moreover, we can assume
$\varepsilon$ small enough to have $\varepsilon CT\leq\eta
e^{-CT}/4$. So, on the event $D^{\varepsilon,y}$, by Gronwall's Lemma,
for $t\geq t_{\psi}-\rho$,
\begin{equation*}
  \|{X}^{\varepsilon,y}_t-\psi_t\|\leq\left(\frac{3}{4}\eta
    e^{-CT}+\sqrt{\varepsilon}\left\|\int_{t_{\psi}-\rho}^t
      \sigma({X}^{\varepsilon,y}_s)dW^\varepsilon_s\right\|\right) e^{CT}.
\end{equation*}
Since $\|{X}^{\varepsilon,y}-\psi\|_{t_{\psi}-\rho,T}\geq\eta$
on $D^{\varepsilon,y}$, we finally obtain
\begin{equation*}
  D^{\varepsilon,y}\subset\left\{\left\|\sqrt{\varepsilon}
      \int_{t_{\psi}-\rho}^{\cdot}\sigma({X}^{\varepsilon,y}_s)dW^\varepsilon_s
    \right\|_{t_{\psi}-\rho,T}\geq\frac{\eta e^{-CT}}{4},\ 
    \|\sqrt{\varepsilon}(W^\varepsilon_{\cdot}-W^\varepsilon_{t_{\psi}-\rho})\|_{t_{\psi}-
      \rho,T}\leq 2\delta\right\}
\end{equation*}
Equation~(\ref{eq:lem-phi=0-(i)}) is now a consequence of
Lemma~\ref{lem:prelim}.\hfill$\Box$

\paragraph{Proof of Lemma~\ref{lem:phi=0}~(iii)}
As for Point~(i), take $\alpha>0$ such that
$S(\varphi)_t\in\Gamma_{\alpha}$ for all $t\in[0,T]$. Fix
$\eta\leq\alpha/2$. Then, on the event
$\{\|{X}^{\varepsilon,y}-S(\varphi)\|_{0,T}\leq\eta\}$, for any
$t\in[0,T]$, ${X}^{\varepsilon,y}_t\in \Gamma_{\alpha/2}$.  Take $C$
such that $b$ and $\sigma$ are $C$-Lipschitz and $\tilde{b}$ is
bounded by $C$ on $\Gamma_{\alpha/2}$. It follows
from~(\ref{eq:Xtilde}) that, on the event
$\{\|{X}^{\varepsilon,y}-S(\varphi)\|_{0,T} \leq\eta\}$, for any
$t\in[0,T]$,
\begin{align*}
  \sqrt{\varepsilon}\biggl\|\int_0^t
  \sigma({X}^{\varepsilon,y}_s)dW^\varepsilon_s\biggr\| & \leq
    \|{X}^{\varepsilon,y}_t-S(\varphi)_t\|+\|y-x\|+\left\|\int_0^t
    [\sigma({X}^{\varepsilon,y}_s)-\sigma(S(\varphi)_s)]\dot
    {\varphi}_sds\right\| \\ &
  \qquad +\left\|\int_0^t[b({X}^{\varepsilon,y}_s)-
    b(S(\varphi)_s)]ds\right\|-\left\|\varepsilon
    \int_0^t\tilde{b}({X}^{\varepsilon,y}_s)ds\right\| \\ &
  \leq 2\eta
  +C\int_0^T(1+\|\dot{\varphi}_s\|)
  \|{X}^{\varepsilon}_s-S(\varphi)_s\|ds+\varepsilon CT \\
  & \leq\eta(2+2CT+C\sqrt{uT})
\end{align*}
if $\varepsilon<\eta$. Therefore,
\begin{multline}
  \{\|{X}^{\varepsilon,y}-S(\varphi)\|_{0,T}\leq\eta,\ 
  \|\sqrt{\varepsilon}W^\varepsilon\|_{0,T}\geq\delta\} \\ \subset
  \left\{\forall t\in[0,T],\ {X}^{\varepsilon,y}_t\in
    \Gamma_{\frac{\alpha}{2}},\ \sqrt{\varepsilon}\left\|\int_0^t
      \sigma({X}^{\varepsilon,y}_s)dW^\varepsilon_s\right\|_{0,T}\leq\eta\beta,\ 
    \sqrt{\varepsilon}\|W^\varepsilon\|_{0,T}\geq\delta\right\}, \label{eq:GL(ii)1}
\end{multline}
where $\beta=2+2CT+C\sqrt{uT}$.

Let
\begin{gather*}
  \tau^{\varepsilon,y}=\inf\{ t:d({X}^{\varepsilon,y}_t,\Gamma)
  \leq\alpha/2\}\wedge T, \\
  Y^{\varepsilon,y}_t=\int_0^t\sigma({X}^{\varepsilon,y}_s)dW^\varepsilon_s, \\
  \xi=\chi\sigma^{-1},
\end{gather*}
where $\chi$ is a Lipschitz function from $\RR^d$ to $[0,1]$ such that
$\chi(x)=0$ if $d(x,\Gamma)\leq\alpha/4$ and $\chi(x)=1$ if
$d(x,\Gamma)\geq\alpha/2$. With these notations,~(\ref{eq:GL(ii)1})
implies
\begin{multline*}
  \{\|{X}^{\varepsilon,y}-S(\varphi)\|_{0,T}\leq\eta,\ 
  \|\sqrt{\varepsilon}W^\varepsilon\|_{0,T}\geq\delta\} \\ \subset
  \left\{\sqrt{\varepsilon}\|Y^{\varepsilon,y}\|_{0,T}\leq\eta\beta,\ 
    \sqrt{\varepsilon}\left\|\int_0^{t\wedge\tau^{\varepsilon,y}}
      \xi({X}^{\varepsilon,y}_s)dY^{\varepsilon,y}_s\right\|_{0,T}
    \geq\delta\right\}.
\end{multline*}
Equation~(\ref{eq:lem-phi=0-(ii)}) is now a direct consequence of
Lemma~\ref{lem:prelim}: $\xi$ is Lipschitz and bounded on $\RR^d$ by
Proposition~\ref{prop:a-b-btilde}~(iii), and for any
$t\leq\tau^{\varepsilon,y}$, $\langle
Y^{\varepsilon,y}\rangle_t=\int_0^t a({X}^{\varepsilon,y}_s)ds$ which
is bounded by a constant $A$ independent of $y$ and $\varepsilon$, by
Proposition~\ref{prop:a-b-btilde}~(i).\hfill$\Box$
\bigskip

Let us come to the proof of Lemmas~\ref{lem:prelim}. It is adapted
from the proof of Lemma~1.3 of~\cite{doss-priouret-83}, and makes use
of Lemma~\ref{lem:Stroock}.

\paragraph{Proof of Lemma~\ref{lem:prelim}}
We use a discretization technique: for any $p\in\mathbb{N}$, we define
${X}^{\varepsilon,y,p}_t={X}^{\varepsilon,y}_{k2^{-p}}$, where
$k\in\mathbb{N}$ is such that $k\leq t2^{p}<k+1$. Given $\gamma>0$,
$p\geq 1$ and $\delta>0$, we can write
\begin{equation*}
  \left\{\left\|\sqrt{\varepsilon}\int_0^{\cdot\wedge\tau}
      \xi({X}^{\varepsilon,y}_s) dY_s\right\|_{0,T}\geq\eta,\
    \|\sqrt{\varepsilon}Y\|_{0,T}\leq \delta\right\}\subset
  A^{\varepsilon}\cup B^{\varepsilon}\cup C^{\varepsilon},
\end{equation*}
where
\begin{gather*}
  A^{\varepsilon}=\{\|{X}^{\varepsilon,y}-
  {X}^{\varepsilon,y,p}\|_{0,\tau}\geq\gamma\}, \\
  B^{\varepsilon}=\left\{\|{X}^{\varepsilon,y}-
    {X}^{\varepsilon,y,p}\|_{0,\tau}\leq\gamma,\ 
    \left\|\sqrt{\varepsilon}\int_0^{\cdot\wedge\tau}
      [\xi({X}^{\varepsilon,y}_s)-\xi({X}^{\varepsilon,y,p}_s)]
      dY_s\right\|_{0,T}\geq \frac{\eta}{2}\right\} \\
  \mbox{and}\quad C^{\varepsilon}=\left\{\left\|\sqrt{\varepsilon}
      \int_0^{\cdot\wedge\tau}\xi({X}^{\varepsilon,y,p}_s)dY_s
    \right\|_{0,T}\geq\frac{\eta}{2},\ 
    \|\sqrt{\varepsilon}Y\|_{0,T}\leq\delta\right\}.
\end{gather*}
We will choose $\gamma$ such that
$\mathbb{P}^{\varepsilon,y}(B^{\varepsilon})$ is sufficiently small,
next $p\geq 1$ to control
$\mathbb{P}^{\varepsilon,y}(A^{\varepsilon})$, and finally $\delta>0$
such that $C^{\varepsilon}=\emptyset$.

First, we apply Lemma~\ref{lem:Stroock} with
$Z_t=\sqrt{\varepsilon}[\xi({X}^{\varepsilon,y}_t)-\xi(
{X}^{\varepsilon,y,p}_t)]$. Let $M_{\gamma}:=
\sup_{\|x-y\|\leq\gamma}\|\xi(x)-\xi(y)\|$, which is finite since
$\xi$ is uniformly continuous. Then, on $B^{\varepsilon}$,
$\|Z_t\|\leq\sqrt{\varepsilon} M_{\gamma}$ for all
$t\leq\tau$. Therefore,
\begin{equation*}
  \mathbb{P}^{\varepsilon,y}(B^{\varepsilon})\leq 2d\exp\left(
    -\frac{\eta^2/4}{2dTA\varepsilon M_{\gamma}^2}\right).
\end{equation*}
Now, $M_{\gamma}\rightarrow 0$ when $\gamma\rightarrow 0$ since $\xi$
is absolutely continuous. Therefore, choosing $\gamma$ small enough,
$\varepsilon\ln\mathbb{P}^{\varepsilon,y}(B^{\varepsilon})\leq -2R$ for all
$\varepsilon\leq 1$.

Second, $\gamma>0$ being fixed as above,~(\ref{eq:Xtilde}) yields
\begin{align*}
  & {\mathbb{P}^{\varepsilon,y}(\|{X}^{\varepsilon,y}-{X}^{\varepsilon,y,p}
  \|_{0,\tau}\geq\gamma)} \\ & \qquad
  \leq\sum_{k=0}^{T2^p-1}\mathbb{P}^{\varepsilon,y}\left(\sup_{k2^{-p}\leq
      t\leq(k+1)2^{-p}}\left\|\int_{k2^{-p}\wedge\tau}^{t\wedge\tau}
      \sqrt{\varepsilon}\sigma({X}^{\varepsilon,y}_s)dW^\varepsilon_s\right\|
    \geq\frac{\gamma}{2}\right) \\
  & \qquad\quad +\sum_{k=0}^{T2^p-1}\mathbb{P}^{\varepsilon,y}\left(\sup_{k2^{-p}\leq
      t\leq(k+1)2^{-p}}\left\|\int_{k2^{-p}\wedge\tau}^{t\wedge\tau}
      \left[b^{\varepsilon}({X}^{\varepsilon,y}_s)+
        \sigma({X}^{\varepsilon,y}_s)\dot{\varphi}_s\right]ds
    \right\|\geq\frac{\gamma}{2}\right) \\
  & \qquad\leq\sum_{k=0}^{T2^p-1}\mathbb{P}^{\varepsilon,y}\left(\sup_{k2^{-p}\leq
      t\leq(k+1)2^{-p}}\left\|\int_{k2^{-p}\wedge\tau}^{t\wedge\tau}
      \sqrt{\varepsilon}\sigma({X}^{\varepsilon,y}_s)dW^\varepsilon_s\right\|
    \geq\frac{\gamma}{2}\right) \\
  & \qquad\quad+\sum_{k=0}^{T2^p-1}\mathbb{P}^{\varepsilon,y}\left(C2^{-p}+
    C2^{-p/2}\sqrt{u}\geq\frac{\gamma}{2}\right),
\end{align*}
where $C$ is a bound for $b^{\varepsilon}$ and $\sigma$ and
$u=\int_0^T\|\dot{\varphi}_s\|^2ds<+\infty$. For $p$ big enough, the
second sum of the right-hand side equals $0$. For the first sum,
Lemma~\ref{lem:Stroock} with $\tau=T=2^{-p}$, $Y=W^\varepsilon$, $A=1$,
$R=\gamma/2$ and $B=\sqrt{\varepsilon}C$ gives that
\begin{equation*}
  \mathbb{P}^{\varepsilon,y}\left(\sup_{k2^{-p}\leq
      t\leq(k+1)2^{-p}}\left\|\int_{k2^{-p}\wedge\tau}^{t\wedge\tau}
      \sqrt{\varepsilon}\sigma({X}^{\varepsilon,y}_s)dW^\varepsilon_s\right\|
    \geq\frac{\gamma}{2}\right)\leq 2d\exp\left(-\frac{\gamma^2/4}
      {2d2^{-p}C^2\varepsilon}\right)
\end{equation*}
for all $0\leq k<T2^p$.  Therefore, taking $p$ large enough,
$\varepsilon\ln\mathbb{P}^{\varepsilon,y} (A^{\varepsilon})\leq -2R$
for all $\varepsilon\leq 1$.

Third, with $p\geq 1$ and $\gamma>0$ as above, for $t\leq T$,
\begin{equation*}
  \sqrt{\varepsilon}\int_0^{t\wedge\tau}\xi({X}^{
    \varepsilon,y,p}_s)dY_s=\sum_{i=0}^{T2^p-1}\sqrt{\varepsilon}
  \xi({X}^{\varepsilon,y}_{i2^{-p}\wedge\tau})
  [Y_{(i+1)2^{-p}\wedge t\wedge\tau}-Y_{i2^{-p}\wedge t\wedge\tau}].
\end{equation*}
Therefore, since
$\|\sqrt{\varepsilon}Y\|_{[0,T]}\leq\delta$ on the event
$C^{\varepsilon}$, we have for all $t\leq T$
\begin{equation*}
  \left\|\sqrt{\varepsilon}\int_0^{t\wedge\tau}\xi
    ({X}^{\varepsilon,y,p}_s)dY_s\right\|
  \leq\sum_{i=0}^{T2^p-1}2\delta C,
\end{equation*}
where $C$ is a bound for $\xi$. Hence $C^{\varepsilon}=\emptyset$ as
soon as $\delta<\eta2^{-(p+2)}/CT$.

We finally obtain that
$\varepsilon\ln\mathbb{P}^{\varepsilon,y}(A^{\varepsilon}\cup
B^{\varepsilon}\cup C^{\varepsilon})\leq\varepsilon\ln 2-2R$, which
yields~(\ref{eq:lem-prelim}) for $\varepsilon$ small enough.

This argument is true for any $y\in\RR^d$ and for any
stopping time $\tau$. It remains to observe that $A$ is the only
information about $Y$ that we used to estimate
$\mathbb{P}^{\varepsilon,y}(B^{\varepsilon})$, that $Y$ does not appear in
$A^{\varepsilon}$, and that no assumption about $Y$ is necessary to
obtain $C^{\varepsilon}=\emptyset$. Hence, the constant $A$ is the
only information about $Y$ required to obtain $\delta$ and
$\varepsilon_0$.\hfill$\Box$
\bigskip

The proof of Lemma~\ref{lem:approx} is now completed. It only remains
to prove Lemmas~\ref{lem:S-I-Itilde}.

\paragraph{Proof of Lemma~\ref{lem:S-I-Itilde}}
Let us first prove Point~(i). Take $\psi\in\tilde{\cal
  C}^{ac}_x([0,T],\RR^d)$. Any $\varphi\in{\cal
  C}^{ac}_0([0,T],\mathbb{R}^d)$ such that $S(\varphi)=\psi$ must
satisfy for any $t\in[0,t_\psi)$
\begin{equation*}
  \dot{\psi}_t=b(\psi_t)+\sigma(\psi_t)\dot{\varphi}_t.
\end{equation*}
Therefore, such a $\varphi$ is uniquely defined for $t<t_{\psi}$ by
\begin{equation}
  \label{eq:I=infJ-1}
  \dot{\varphi}_t=\sigma^{-1}(\psi_t)[\dot{\psi}_t-b(\psi_t)].
\end{equation}
Thus
\begin{equation*}
  I_{T,x}(\psi)=\frac{1}{2}\int_0^{t_{\psi}}\|\sigma^{-1}(\psi_t)
  [\dot{\psi}_t-b(\psi_t)]\|^2dt=\frac{1}{2}\int_0^{t_{\psi}}\|\dot{\varphi}_t\|^2dt
  \leq J_T(\varphi)
\end{equation*}
for any $\varphi$ such that $S(\varphi)=\psi$, and
$I_{T,x}(\psi)=J_T(\varphi)$ if and only if $\dot{\varphi}_t=0$ for
all $t>t_\psi$.

This trivially implies that $I_{T,x}(\psi)=\inf\{J_T(\varphi),\:
S(\varphi)=\psi\}$ when $I_{T,x}(\psi)=+\infty$. In the case where
$I_{T,x}(\psi)<+\infty$, we clearly have $I_{T,x}(\psi)\leq
\inf\{J_T(\varphi),\: S(\varphi)=\psi\}$. To prove the converse
inequality, it suffices to check that there exists an absolutely
continuous function $\varphi$ satisfying~(\ref{eq:I=infJ-1}) for
$t<t_\psi$ and $\dot{\varphi}_t=0$ for $t\geq t_\psi$. This is
equivalent to the fact that
$\sigma^{-1}(\psi_t)[\dot{\psi}_t-b(\psi_t)]$ is ${L}^1$ on
$[0,t_{\psi}]$. Since $I_{T,x}(\psi)<+\infty$, this function is
actually $L^2$, which ends the proof of Point~(i).

For Point~(ii), remind that $\sigma$ is uniformly non-degenerate on
$\Gamma_{\alpha}$ for any $\alpha>0$. Therefore, the fact that ${\cal
  C}^1([0,T],\RR^d\setminus\Gamma)\subset S(\{J_T<\infty\})$ follows
from~(\ref{eq:I=infJ-1}). Since $S(\{J_T<\infty\})\subset \tilde{\cal
  C}^{ac}_x([0,T],\RR^d)$ and any function of $\tilde{\cal
  C}^{ac}_x([0,T],\RR^d)$ is the limit of elements of ${\cal
  C}^1([0,T],\RR^d\setminus\Gamma)$, Point~(ii) is clear.\hfill$\Box$

\section{Application to the problem of exit from a domain}
\label{sec:exitdomain}

We study in this section the biological phenomenon of punctualism. We
consider a bounded open subset $G$ of $\RR^d$ containing a unique,
stable equilibrium of the canonical equation of adaptive dynamics
$\dot\phi=b(\phi)$. We will assume for convenience that this
equilibrium is 0. Note that the equilibria of the canonical equation
are exactly the points of $\Gamma$. As observed in
Remark~\ref{rem:canonical-eq}, when $\varepsilon$ is small,
$X^{\varepsilon,x}$ is close to the solution of this ODE with initial
state $x$ with high probability. Yet, the diffusion phenomenon may
almost surely drive $X^{\varepsilon,x}$ out of $G$. Our next result
gives estimates of the time and position of exit of $X^{\varepsilon}$
from $G$ (``problem of exit from a
domain''~\cite{freidlin-wentzell-84}).

We will follow closely section~5.7 of Dembo and
Zeitouni~\cite{dembo-zeitouni-93}, where a similar result for
non-degenerate diffusions is proved.

When the initial condition of the solution of the
SDE~(\ref{eq:SDE-LDP}) constructed in Proposition~\ref{prop:exist-LDP}
is not precised, it will by denoted by $X^\varepsilon$. The value of
$X^\varepsilon$ at time 0 will then be specified by considering the
probability of events involving $X^\varepsilon$ under $\PP_x$, which
is the law of the process $X^{\varepsilon,x}$. Expectations with
respect to $\PP_x$ will be denoted $\EE_x$. We will also use
throughout this section the notation
$B(\rho):=\{y\in\mathbb{R}^d:\|y\|\leq\rho\}$ and
$S(\rho)=\{y\in\mathbb{R}^d:\|y\|=\rho\}$ for $\rho>0$. It will always
be implicitly assumed that $\rho>0$ is small enough to have
$B(\rho)\subset G$ and $S(\rho)\subset G$.

We will assume $d\geq 2$. Otherwise, the problem has few interest: if
$G=(c,c')\subset\RR$ contains a unique point $x$ of $\Gamma$, and if
$y>x$ (say), the process $X^{\varepsilon,y}$ can exit $G$, only at
$c'$, and the probability of reaching $x$ before $c'$ can be computed
explicitly using classical results on one-dimensional diffusion
processes~\cite[Prop.\:5.5.22]{karatzas-shreve-88}.
\medskip

Let
\begin{equation*}
  V(y,z,t)=\inf_{\{\psi\in{\cal C}([0,t],\RR^d):\psi(0)=y,
    \psi(t)=z\}}\tilde{I}_{t,y}(\psi),
\end{equation*}
which is, heuristically, the cost of forcing $X^{\varepsilon,y}$ to be
at $z$ at time $t$. Define also
\begin{equation*}
  V(y,z)=\inf_{t>0}\:V(y,z,t).
\end{equation*}
The function $V(0,z)$ is called the
\emph{quasi-potential}~\cite{freidlin-wentzell-84}.

Six assumptions are required for our result:
\begin{description}
\item[\textmd{(Ha)}] $G$ is a bounded open subset of $\RR^d$ such that
  $G\cap\Gamma=\{0\}$ and with sufficiently smooth boundary $\partial
  G$ for
  \begin{equation*}
    \tau^{\varepsilon}=\inf\{t>0:X^{\varepsilon}_t\in\partial G\}
  \end{equation*}
  to be a well-defined stopping time. Moreover, for any solution of
  \begin{equation}
    \label{eq:deterministic}
    \dot{\phi}=b(\phi)
  \end{equation}
  such that $\phi(0)\in G$, we have $\phi(t)\in G$ for all $t>0$ and
  $\lim_{t\rightarrow\infty}\phi(t)=0$.
\item[\textmd{(Hb)}] $\bar{V}:=\inf_{z\in\partial G}{V}(0,z)<\infty$.
\item[\textmd{(Hc)}] For any $\varepsilon>0$ and $y\in
  G\setminus\{0\}$, $\displaystyle{\mathbb{P}_y\left(
      \lim_{t\rightarrow\infty}X^{\varepsilon}_t=0\right)=0}$.
\item[\textmd{(Hd)}] The points of $\Gamma$ are isolated in $\RR^d$.
\item[\textmd{(He)}] For any $y\in\overline{G}\cap\Gamma$, $g$ is
  ${\cal C}^2$ at $(y,y)$ and $H_{1,1}g(y,y)+H_{1,2}g(y,y)$ is
  invertible.
\item[\textmd{(Hf)}] All the trajectories of the deterministic
  system~(\ref{eq:deterministic}) with initial value in $\partial G$
  converges to $0$ as $t\rightarrow\infty$.
\end{description}

Assumption~(Ha) states that the domain $G$ is an \emph{attracting}
domain for~(\ref{eq:deterministic}). If Assumption~(Hb) fails, all
points of $\partial G$ are equally unlikely on the large deviations
scale. We have given in Theorem~\ref{thm:d>1}
(sections~\ref{sec:cond-tau-infinite}) conditions under which~(Hc)
holds. Assumption~(Hd) is required in the large deviation
Theorem~\ref{thm:true-LDP}. We have already encountered an assumption
similar to~(He) in Propositions~\ref{prop:d>1}
and~\ref{prop:I-not-lsc}. It allows to control the non-degeneracy of
$a(x)$ near $\overline{G}\cap\Gamma$. Finally, Assumption~(Hf)
prevents situations where $\partial G$ is the characteristic boundary
of the domain of attraction of $0$. This last assumption is needed
only for Point~(b) of the following result. Note that when (Hf) is
true, $\overline{G}\cap\Gamma=\{0\}$

\begin{thmbr}
  \label{thm:exit-domain}
  \begin{description}
  \item[\textmd{(a)}] Assume (H) and (Ha--e). Then, for all
    $x\in G\setminus\{0\}$ and $\delta>0$,
    \begin{equation}
      \label{eq:thm-exit-domain-(a)-1}
      \lim_{\varepsilon\rightarrow 0}\mathbb{P}_x(
        \tau^{\varepsilon}>e^{(\bar{V}
        -\delta)/\varepsilon})=1.
    \end{equation}
  \item[\textmd{(b)}] Assume (H) and (Ha--f).  If $N$ is a
    closed subset of $\partial G$ and if $\inf_{z\in
      N}{V}(0,z)>\bar{V}$, then for any $x\in G\setminus\{0\}$,
    \begin{equation}
      \label{eq:thm-exit-domain-(b)-1}
      \lim_{\varepsilon\rightarrow
        0}\mathbb{P}_x(X^{\varepsilon}_{\tau^{\varepsilon}}\in N)=0.
    \end{equation}
    In particular, if there exists $z^*\in\partial G$ such that
    ${V}(0,z^*)<{V}(0,z)$ for all $z\in\partial
    G\setminus\{z^*\}$, then, for any $\delta>0$ and $x\in G\setminus\{0\}$,
    \begin{equation}
      \label{eq:thm-exit-domain-(b)-2}
      \lim_{\varepsilon\rightarrow
      0}\mathbb{P}_x(\|X^{\varepsilon}_{\tau^{\varepsilon}}-z^*\|<\delta)=1.
    \end{equation}
  \end{description}
\end{thmbr}

The proof of such results is classically guided by the heuristics that, as
$\varepsilon\rightarrow 0$, $X^{\varepsilon}$ wanders around $0$ for
an exponentially long time, during which its chance of hitting a
closed set $N\subset\partial G$ is determined by $\inf_{z\in
  N}{V}(0,z)$. Any excursion off the stable point $0$ has an
overwhelmingly high probability of being pulled back near $0$, and it
is not the time spent near any part of $\partial G$ that matters but
the \emph{a priori} chance for a direct, fast exit due to a rare
segment in the Brownian motion's path.

Usually, such results also include an upper bound for
$\tau^\varepsilon$. We are not able to obtain such a result because of
the singularity of the process $X^\varepsilon$ at 0. Because the
matrix $a(x)$ is 0 at $x=0$, the time spent by the process near 0
before hitting $S(\rho)$ is not uniformly bounded (in probability)
with respect to the initial condition (actually, it is even infinite
when $X^\varepsilon_0=0$).

For this reason, the proof of a similar result in Dembo and
Zeitouni~\cite{dembo-zeitouni-93} (Thm.\:5.7.11 and Cor.\:5.7.16)
cannot be directly adapted to our situation. Below, we are only going
to detail the steps that must be modified.  In particular,
Theorem~\ref{thm:exit-domain}~(a) will be obtained exactly as
in~\cite{dembo-zeitouni-93}, whereas Point~(b) has to be obtained
without using any upper bound on $\tau^{\varepsilon}$. 

We are going to use four lemmas. The first one gives estimates on
continuity of $V(x,\cdot,t)$ around 0 and $\partial G$.
\begin{lemma}
  \label{lem:V-cont}
  Assume (H),~(Hd) and~(He). For any $\delta>0$, there exists $\rho>0$
  small enough such that
  \begin{equation}
    \label{eq:lem-V-cont-1}
    \sup_{(x,y)\in (B(\rho)\setminus\{0\})\times B(\rho)}\:\inf_{t\in[0,1]}V(x,y,t)<\delta
  \end{equation}
  and
  \begin{equation}
    \label{eq:lem-V-cont-2}
    \sup_{\{(x,y)\in(\RR^d\setminus\Gamma)\times\RR^d,\ \inf_{z\in\partial
    G}(\|y-z\|+\|x-z\|)\leq\rho\}}\:\inf_{t\in[0,1]}V(x,y,t)<\delta.
  \end{equation}
\end{lemma}
For the next lemmas, we define
\begin{equation*}
  \sigma_{\rho}:=\inf\{t\geq 0:X^{\varepsilon}\in
  B(\rho)\cup\partial G\}.
\end{equation*}
The second lemma gives a uniform lower bound on the probability of an
exit from $G$ starting from a small sphere around $0$ before hitting
an even smaller sphere.
\begin{lemma}
  \label{lem:LB-exit-G}
  Assume~(H) and (Ha--e). Then
  \begin{equation*}
    \lim_{\rho\rightarrow 0}\liminf_{\varepsilon\rightarrow 0}
    \varepsilon\ln\inf_{y\in S(2\rho)}\mathbb{P}_y(X^{\varepsilon}_{
      \sigma_{\rho}}\in\partial G)\geq -\bar{V}.
  \end{equation*}
\end{lemma}
The following upper bound relates the quasi-potential
${V}(0,\cdot)$ with the probability that an excursion starting
from a small sphere around $0$ hits a given subset of $\partial G$
before hitting an even smaller sphere.
\begin{lemma}
  \label{lem:UB-exit-N}
  Assume~(H) and (Ha--f). For any closed set
  $N\subset\partial G$,
  \begin{equation*}
    \lim_{\rho\rightarrow 0}\limsup_{\varepsilon\rightarrow 0}
    \varepsilon\ln\sup_{y\in S(2\rho)}\mathbb{P}_y(X^{\varepsilon}_{
      \sigma_{\rho}}\in N)\leq -\inf_{z\in N}{V}(0,z)
  \end{equation*}
\end{lemma}
The last lemma is used to extend the previous upper bound to
any initial condition $x\in G$.
\begin{lemma}
  \label{lem:return-to-0}
  Assume~(H) and (Ha). For every $\rho>0$ such that
  $B(\rho)\subset G$ and all $x\in G$,
  \begin{equation*}
    \lim_{\varepsilon\rightarrow 0}\mathbb{P}_x(X^{\varepsilon}_{
      \sigma_{\rho}}\in B(\rho))=1.
  \end{equation*}
\end{lemma}

The statements of Lemmas~\ref{lem:V-cont},~\ref{lem:UB-exit-N}
and~\ref{lem:return-to-0} are the same as Lemmas~5.7.8,~5.7.21
and~5.7.22 of~\cite{dembo-zeitouni-93}, respectively. Among them,
Lemmas~\ref{lem:UB-exit-N} and~\ref{lem:return-to-0} can be deduced
from Corollary~\ref{cor:LDP} exactly as in~\cite{dembo-zeitouni-93},
so we omit their proof. Because of the degeneracy of $X^\varepsilon$
at 0, Lemma~\ref{lem:V-cont} must be proved with a different
method. Finally, Lemma~\ref{lem:LB-exit-G} replaces Lemma~5.7.18
of~\cite{dembo-zeitouni-93} and is very different since it gives no
upper control on $\tau^\varepsilon$. This lemma and the proof of
Theorem~\ref{thm:exit-domain}~(b) are the new part of our proof.

Theorem~\ref{thm:exit-domain}~(a) can be proved exactly as the
corresponding inequalities in Theorem~5.7.11 and Corollary~5.7.16
of~\cite{dembo-zeitouni-93}. It makes use of our
Lemmas~\ref{lem:V-cont},~\ref{lem:UB-exit-N}
and~\ref{lem:return-to-0}, and of Lemmas~5.7.19 and~5.7.23
of~\cite{dembo-zeitouni-93}, which can be proved exactly as therein.
One simply must take care that $x$ belongs to $G\setminus\{0\}$
instead of $G$.  Let us omit this proof.

We first give the proof of Theorem~\ref{thm:exit-domain}~(b) and next
those of Lemmas~\ref{lem:V-cont} and~\ref{lem:LB-exit-G}.

\paragraph{Proof of Theorem~\ref{thm:exit-domain}~(b)}
Let $\rho>0$ be small enough to have $B(2\rho)\subset G$ (the precise
choice of $\rho$ will be specified later). Let $\theta_0=0$ and for
$m=0,1,\ldots$ define the stopping times
\begin{equation}
  \label{eq:def-tau-m-theta-m}
  \begin{aligned}
    \tau_m & =\inf\{t\geq\theta_m:X^{\varepsilon}_t\in
    B(\rho)\cup\partial G\}, \\ \theta_{m+1} &
    =\inf\{t>\tau_m:X^{\varepsilon}_t\in S(2\rho)\},
  \end{aligned}
\end{equation}
with the convention that $\theta_{m+1}=\infty$ if
$X^{\varepsilon}_{\tau_m}\in\partial G$. Each interval
$[\tau_m,\tau_{m+1}]$ represents one significant excursion off
$B(\rho)$. Note that, necessarily, $\tau^{\varepsilon}=\tau_m$ for
some integer $m$.

First, Assumption~(Hc) implies that $\theta_{m+1}<\infty$ as soon as
$X^{\varepsilon}_{\tau_m}\in B(\rho)$. This can be proved as follows.

On the one hand, Assumption~(Hc) implies that, for all $x\in S(\rho)$,
\begin{equation}
  \label{eq:zut-zut}
  \lim_{\alpha\rightarrow 0}\mathbb{P}_x(\limsup_{t\rightarrow
    +\infty}\|X^{\varepsilon}_t\|\geq\alpha)=1.
\end{equation}
On the other hand, for any $\alpha>0$, $X^{\varepsilon}$ is a
diffusion with bounded drift part and uniformly non-degenerate
diffusion part in $B(2\rho)\cap\Gamma_{\alpha/2}$. Therefore,
$X^{\varepsilon}$ has a uniformly positive probability to reach
$S(2\rho)$ before $S(\alpha/2)$ starting from any point of
$S(\alpha)$. Hence, by the strong Markov property of
Proposition~\ref{prop:str-M-pty}, for all $x\in S(\rho)$,
$$
\mathbb{P}_x(\theta_1<\infty\mid\limsup_{t\rightarrow+\infty}
\|X^{\varepsilon}_t\|\geq\alpha)=1.
$$
Combining this with~(\ref{eq:zut-zut}) we have that
$\mathbb{P}_x(\theta_1<\infty)=1$ for all $x\in S(\rho)$, which
implies the required result.

Second, fix a closed set $N\subset G$ such that $\bar{V}_N:=\inf_{z\in
  N}{V}(0,z)>\bar{V}$. Assume $\bar{V}_N<\infty$ (otherwise,
$\bar{V}_N$ may be replaced by any arbitrary large constant in the
proof below). Fix $\eta>0$ such that
$\eta<(\bar{V}_N-\bar{V})/3$. Applying Lemmas~\ref{lem:LB-exit-G}
and~\ref{lem:UB-exit-N}, we fix $\rho>0$ and $\varepsilon_0>0$ such
that
\begin{equation}
  \label{eq:pf-thm-exit-(b)-0}
  \inf_{y\in S(2\rho)}\mathbb{P}_y(X^{\varepsilon}_{\sigma_{\rho}}
  \in\partial G)\geq e^{-(\bar{V}+\eta)/\varepsilon},
  \quad\forall\varepsilon\leq\varepsilon_0
\end{equation}
and
\begin{equation*}
  \sup_{y\in S(2\rho)}\mathbb{P}_y(X^{\varepsilon}_{\sigma_{\rho}}\in
  N)\leq e^{-(\bar{V}_N-\eta)/\varepsilon},\quad\forall\varepsilon
  \leq\varepsilon_0.
\end{equation*}
Fix $y\in B(\rho)$. For any $l\geq 1$, we have
\begin{equation}
  \label{eq:pf-thm-exit-(b)-1}
  \mathbb{P}_y(X^{\varepsilon}_{\tau^{\varepsilon}}\in
  N)\leq\mathbb{P}_y(\tau^{\varepsilon}>\tau_l)
  +\sum_{m=1}^l\mathbb{P}_y(\tau^{\varepsilon}=\tau_m\mbox{\ and\
    }X^{\varepsilon}_{\tau^{\varepsilon}}\in N).
\end{equation}
The second term can be bounded as follows: for $m\geq 1$, $y\in
B(\rho)$ and $\varepsilon\leq\varepsilon_0$, it follows from the
strong Markov property that
\begin{align*}
  \mathbb{P}_y(\tau^{\varepsilon}=\tau_m\mbox{\ and\ 
  }X^{\varepsilon}_{\tau^{\varepsilon}}\in N) &
  =\mathbb{P}_y(\tau^{\varepsilon}>\tau_{m-1})\mathbb{P}_y
  (X^{\varepsilon}_{\tau_m}\in N\mid\tau^{\varepsilon}>\tau_{m-1}) \\
  & =\mathbb{P}_y(\tau^{\varepsilon}>\tau_{m-1})\mathbb{E}_y
  [\mathbb{P}_{X^{\varepsilon}_{\theta_m}}(X^{\varepsilon}_{\sigma_{\rho}}\in
  N)\mid\tau^{\varepsilon}>\tau_{m-1}] \\
  & \leq\sup_{x\in
    S(2\rho)}\mathbb{P}_x(X^{\varepsilon}_{\sigma_{\rho}}\in N)\leq
  e^{-(\bar{V}_N-\eta)/\varepsilon}.
\end{align*}
Concerning the first term of the right-hand side
of~(\ref{eq:pf-thm-exit-(b)-1}), for any $l\geq 1$ and $y\in B(\rho)$,
\begin{equation*}
  \mathbb{P}_y(\tau^{\varepsilon}>\tau_l)=\mathbb{E}_y
  [\mathbb{P}_{X^{\varepsilon}_{\theta_1}}(\tau^{\varepsilon}>\tau_{l-1})]
  \leq\sup_{x\in S(2\rho)}\mathbb{P}_x(\tau^{\varepsilon}>\tau_{l-1}),
\end{equation*}
and, for any $x\in S(2\rho)$ and $k\geq 1$,
\begin{align*}
  \mathbb{P}_x(\tau^{\varepsilon}>\tau_k) &
  =[1-\mathbb{P}_x(\tau^{\varepsilon}=\tau_k
  \mid\tau^{\varepsilon}>\tau_{k-1})]\mathbb{P}_x(\tau^{\varepsilon}>\tau_{k-1})
  \\ & =[1-\mathbb{E}_x[\mathbb{P}_{X^{\varepsilon}_{\theta_k}}
  (X^{\varepsilon}_{\sigma_{\rho}}\in\partial G)
  \mid\tau^{\varepsilon}>\tau_{k-1}]]\mathbb{P}_x(\tau^{\varepsilon}>\tau_{k-1})
  \\ & \leq(1-q)\mathbb{P}_x(\tau^{\varepsilon}>\tau_{k-1}),
\end{align*}
where $q:=\inf_{y\in S(2\rho)}\mathbb{P}_y
(X^{\varepsilon}_{\sigma_{\rho}}\in\partial G)\geq
e^{-(\bar{V}+\eta)/\varepsilon}$ by~(\ref{eq:pf-thm-exit-(b)-0}).
Therefore,
\begin{equation*}
  \sup_{y\in S(2\rho)}\mathbb{P}_x(\tau^{\varepsilon}>\tau_k)\leq(1-q)^k.
\end{equation*}
Putting together these estimates in~(\ref{eq:pf-thm-exit-(b)-1}), we
obtain that, for all $y\in B(\rho)$ and $\varepsilon\leq\varepsilon_0$
\begin{equation*}
  \mathbb{P}_y(X^{\varepsilon}_{\tau^{\varepsilon}}\in N)\leq
  \left(1-e^{-\frac{\bar{V}+\eta}{\varepsilon}}\right)^{l-1}
  +le^{-\frac{\bar{V}_N-\eta}{\varepsilon}}.
\end{equation*}
We choose $l=\lfloor 2e^{(\bar{V}+2\eta)/\varepsilon}\rfloor$, where
$\lfloor\cdot\rfloor$ denotes the integer part function. Then, for
$\varepsilon$ small enough, $l-1>e^{(\bar{V}+2\eta)/\varepsilon}$ and
\begin{equation*}
  \mathbb{P}_y(X^{\varepsilon}_{\tau^{\varepsilon}}\in N)\leq
  \left(\left(1-\frac{1}{u_{\varepsilon}}\right)^{u_{\varepsilon}}
  \right)^{e^{\eta/\varepsilon}}
  +2e^{\frac{\bar{V}-\bar{V}_N+3\eta}{\varepsilon}},
\end{equation*}
where $u_{\varepsilon}:=e^{(\bar{V}+\eta)/\varepsilon}$. Since
$u_\varepsilon\rightarrow +\infty$, we have
$(1-1/u_{\varepsilon})^{u_{\varepsilon}}\rightarrow 1/e$, and,
finally,
$$
\lim_{\varepsilon\rightarrow 0}\:\sup_{y\in B(\rho)}
\mathbb{P}_y(X^{\varepsilon}_{\tau^{\varepsilon}}\in
N)= 0.
$$
The proof of~(\ref{eq:thm-exit-domain-(b)-1}) is now completed by
combining Lemma~\ref{lem:return-to-0} and the inequality
\begin{equation*}
  \mathbb{P}_x(X^{\varepsilon}_{\tau^{\varepsilon}}\in N)\leq
  \mathbb{P}_x(X^{\varepsilon}_{\sigma_{\rho}}\not\in B(\rho))
  +\sup_{y\in B(\rho)}\mathbb{P}_y
  (X^{\varepsilon}_{\tau^{\varepsilon}}\in N).
\end{equation*}

Applying~(\ref{eq:thm-exit-domain-(b)-1}) to $N=\{z\in\partial
G:\|z-z^*\|\geq\delta\}$ and observing that Lemma~\ref{lem:V-cont}
implies the continuity of $z\mapsto{V}(0,z)$ on $\partial G$, we
easily obtain~(\ref{eq:thm-exit-domain-(b)-2}).\hfill$\Box$

\paragraph{Proof of Lemma~\ref{lem:V-cont}~(\ref{eq:lem-V-cont-1})}
Fix $\delta,\rho>0$, $x\in B(\rho)\setminus\{0\}$ and $y\in
B(\rho)$. In order to simplify the notations, we will use the complex
notation for the coordinates of points of the (two-dimensional) plane
of $\mathbb{R}^d$ containing $0$, $x$ and $y$, and we will assume that
$x=r\in\mathbb{R}$ and $y=r'e^{i\theta}$, with $0<r\leq\rho$ and
$0\leq r'\leq\rho$. Define $\psi\in{\cal C}([0,1],B(\rho))$ by
\begin{equation*}
  \psi(t)=\left\{
    \begin{array}{ll}
      (1-(3t)^2)r+(3t)^2\rho & \mbox{if}\quad 0\leq t\leq 1/3 \\
      \rho e^{i\theta(3t-1)} & \mbox{if}\quad 1/3\leq t\leq 2/3 \\
      (1-(3-3t)^2)r'e^{i\theta}+(3-3t)^2\rho e^{i\theta} &
      \mbox{if}\quad 2/3\leq t\leq 1.
    \end{array}\right.
\end{equation*}
Then $\psi(0)=x$ and $\psi(1)=y$, and $\psi(t)\in
B(\rho)\setminus\{0\}$ for any $t\in[0,1)$. Moreover, for $0\leq t\leq
1/3$, $\psi(t)=r+9t^2(\rho-r)$, so that $\|\psi(t)\|\geq
9t^2(\rho-r)$, and, similarly, for $2/3\leq t\leq 1$, $\|\psi(t)\|\geq
9(1-t)^2(\rho-r')$. Thanks to assumption~(He), a calculation similar
to equation~(\ref{eq:calc-prop-4.1}) in the proof of
Proposition~\ref{prop:I-not-lsc} gives that, with the same $K$, ${\cal
  N}_0$ and $a_0$ as therein, if $B(\rho)\subset{\cal N}_0$,
\begin{align*}
  I_{1,x}(\psi) & \leq\frac{1}{2a_0}\left(\int_0^{1/3}
    \frac{2(18t(\rho-r))^2+2K^2\|\psi(t)\|^2}{\|\psi(t)\|}dt\right. \\
  & \phantom{\leq\frac{1}{2a_0}}\quad+\int_{1/3}^{2/3}\frac{2(3\theta\rho)^2
    +2K^2\|\psi(t)\|^2}{\|\psi(t)\|} \\
  & \left.\phantom{\leq\frac{1}{2a_0}}\quad+\int_{2/3}^1\frac{2(18(1-t)(\rho-r'))^2
      +2K^2\|\psi(t)\|^2}{\|\psi(t)\|}dt\right) \\
  & \leq\frac{1}{2a_0}\left(\int_0^{1/3}
    (648(\rho-r)+2K^2\|\psi(t)\|)dt
    +\int_{1/3}^{2/3}(18\theta^2+2K^2)\rho dt\right. \\
  & \left.\phantom{\leq\frac{1}{2a_0}}\quad+\int_{2/3}^1(648(\rho-r')
    +2K^2\|\psi(t)\|)dt\right) \\
  & \leq\frac{(216+2K^2/3)\rho+(6\theta^2+2K^2/3)\rho
    +(216+2K^2/3)\rho}{2a_0}.
\end{align*}
Consequently, for sufficiently small $\rho>0$ not depending on $x$ and
$y$, $I_{1,x}(\psi)\leq\delta/2$, which
yields~(\ref{eq:lem-V-cont-1}).\hfill$\Box$

\paragraph{Proof of Lemma~\ref{lem:V-cont}~(\ref{eq:lem-V-cont-2})}
Fix $\delta>0$. Thanks to Assumption~(He), using the same method as
above, for any $z\in\partial G\cap\Gamma$, one can find a positive
$\rho_z$ such that
\begin{equation}
  \label{eq:pf-V-cont}
  \sup_{(x,y)\in (B(z,\rho_z)\setminus\{0\})\times
    B(z,\rho_z)}\:\inf_{t\in[0,1]}V(x,y,t)<\delta/2,
\end{equation}
where $B(z,r)$ is the closed ball centered at $z$ with radius $r$

Let $\bar{\rho}_0$ be the infimum of $\rho_z$ for $z\in\partial
G\cap\Gamma$. Since $G$ is bounded, because of Assumption~(Hd), this
set is finite and $\bar{\rho}_0>0$. Reducing $\bar{\rho}_0$ if
necessary, we can assume that $B(\bar{\rho}_0)\subset G$ and that
$d(\Gamma\cap(\RR^d\setminus\overline{G}),\overline{G})>\bar{\rho}_0$.

Fix $x$ and $y$ in $\RR^d\setminus\bigcup_{z\in\partial
  G\cap\Gamma}B(z,\bar{\rho}_0)$ and assume that there exists
$z\in\partial G$ with $\|x-z\|+\|y-z\|\leq\bar{\rho}_0/3$. Then
$d(x,\Gamma)>2\bar{\rho}_0/3$ and
$d(y,\Gamma)>2\bar{\rho}_0/3$. Moreover, since
$\|x-y\|\leq\bar{\rho}_0/3$, the segment $[x,y]$ is included in
$\Gamma_{\bar{\rho}_0/3}$.

Now, for any $t_0>0$, $x$ and $y$ such that
$[x,y]\subset\Gamma_{\bar{\rho}_0/3}$, define $\psi^{(t_0)}\in{\cal
  C}([0,t_0],\RR^d)$ by
\begin{equation*}
  \psi^{(t_0)}(t)=\left(1-\frac{t}{t_0}\right)x+\frac{t}{t_0}y
\end{equation*}
for $0\leq t\leq t_0$. Then $\psi^{(t_0)}(0)=x$ and
$\psi^{(t_0)}(t_0)=y$ and $\psi^{(t_0)}(t)\in\Gamma_{\bar{\rho}_0/3}$
for all $t\in[0,t_0]$.

Since $a$ is uniformly non-degenerate on $\Gamma_{\bar{\rho}_0/3}$,
there exists a constant $C$ bounding the eigenvalues of $a^{-1}$ on
this set. Then
\begin{align*}
  I_{t_0,x}(\psi^{(t_0)}) & \leq\frac{C}{2}\int_0^{t_0}
  (\|\dot{\psi}^{(t_0)}(t)\|^2+\|b(\psi^{(t_0)}(t))\|^2)dt \\
  & \leq \frac{C}{2}\left(\frac{\|x-y\|^2}{t_0}+B^2t_0\right),
\end{align*}
where $B$ is a bound for $b$ on $\RR^d$. Taking $t_0=\|x-y\|/B$, we obtain
\begin{equation*}
  I_{\|x-y\|/B,x}(\psi^{(\|x-y\|/B)})\leq BC\|x-y\|.
\end{equation*}
Therefore, there exists $\bar{\rho}_1>0$ such that
$\inf_{t\in[0,1]}V(x,y,t)<\delta/2$ for any $x$ and $y$ such that
$[x,y]\subset \Gamma_{\bar{\rho}_0/3}$ and
$\|x-y\|\leq\bar{\rho}_1$. In view of~(\ref{eq:pf-V-cont}),
$\rho=\bar{\rho}_1\wedge(\bar{\rho}_0/3)$ is an appropriate constant
in~(\ref{eq:lem-V-cont-2}).\hfill$\Box$

\paragraph{Proof of Lemma~\ref{lem:LB-exit-G}}
Fix $\eta>0$ and let $\rho>0$ be small enough to have $B(2\rho)\subset G$
and for Lemma~\ref{lem:V-cont} to hold with $\delta=\eta/3$ and
$2\rho$ instead of $\rho$. Note that the definition of $\tilde{I}_{t,x}$
implies the inequality $\inf_{y\in S(2\rho)}V(y,z)\leq{V}(0,z)$ as
soon as $z\not\in B(2\rho)$.  

Then, by~(\ref{eq:lem-V-cont-2}) and Assumption~(Hb), there exists
$x\in S(2\rho)$, $z\not\in\overline{G}$, $T_1<\infty$ and
$\psi\in{\cal C}([0,T_1],\RR^d)$ such that $\psi(0)=x$, $\psi(T_1)=z$
and $\tilde{I}_{T_1,x}(\psi)\leq\bar{V}+\eta/3$.  Moreover, by removing the
beginning of the path $\psi$ until the last time where it hits
$S(2\rho)$, we can suppose that for all $t>0$, $\psi(t)\not\in
B(2\rho)$.

Thanks to~(\ref{eq:lem-V-cont-1}), for any $y\in S(2\rho)$, there
exists a continuous path $\psi^y$ of length $t_y\leq 1$ such that
$\psi^y(0)=y$, $\psi^y(t_y)=x$, and
$\tilde{I}_{t_y,y}(\psi^y)\leq\eta/3$.  Moreover, the construction of
this function in the proof of Lemma~\ref{lem:V-cont} allows us to
assume that $\|\psi^y(t)\|=2\rho$ for all $t\in[0,t_y]$. Let $\phi^y$
denote the path obtained by concatenating $\psi^y$ and $\psi$ (in that
order) and extending the resulting function to be of length
$T_0=T_1+1$ by following~(\ref{eq:deterministic}) after reaching
$z$. Since the latter path does not contribute to the rate function,
we obtain that $\tilde{I}_{T_0,y}(\phi^y)\leq\bar{V}+2\eta/3$.

Since $z\in\RR^d\setminus\overline{G}$, the constant
$\Delta:=d(z,\partial G)$ is positive. Define
\begin{equation*}
  O:=\bigcup_{y\in S(2\rho)}\left\{\psi\in{\cal C}([0,T_0],
    \RR^d),\:\|\psi-\phi^y\|_{0,T_0}\leq
    \frac{\Delta\wedge\rho}{2}\right\}.
\end{equation*}
Observe that $O$ is an open subset of ${\cal C}([0,T_0],\RR^d)$ that
contains the functions $\{\phi^y\}_{y\in S(2\rho)}$.  Therefore, by
Corollary~\ref{cor:LDP},
\begin{equation*}
  \liminf_{\varepsilon\rightarrow 0}\varepsilon\ln\inf_{y\in
    S(2\rho)}\mathbb{P}_y(X^{\varepsilon}\in O) \geq -\sup_{y\in
    S(2\rho)}\inf_{\psi\in O}\tilde{I}_{T_0,y}(\psi)
  \geq -\sup_{y\in S(2\rho)}\tilde{I}_{T_0,y}(\phi^y)>-(\bar{V}+\eta).
\end{equation*}
If $\psi\in O$, then $\psi$ reaches the open ball of radius $\Delta/2$
centered at $z$ before hitting $B(\rho)$, so $\psi$ hits $\partial G$
before hitting $B(\rho)$. Hence, for $X^{\varepsilon}_0=y\in
S(2\rho)$, the event $\{X^{\varepsilon}\in O\}$ is contained in
$\{X^{\varepsilon}_{\sigma_{\rho}}\in\partial G\}$, and the proof is
completed.\hfill$\Box$
\bigskip

\textbf{Acknowledgments:} I would like to thank Sylvie M\'el\'eard for
her constant support and many fruitful discussions about this work. I
would also like to thank R\'egis Ferri\`ere for useful discussions
and comments on the biological motivation of this work. Finally, I
would like to thank G\'erard Ben Arous for having pointed out this
problem and for useful comments on a preliminary draft of this work.




\end{document}